\documentclass[11pt,paperA]{article}
\usepackage{srcltx}
\usepackage{amssymb}
\usepackage{amsfonts}
\usepackage{amsmath}
\usepackage{amsthm}
\usepackage{enumerate}
\usepackage{multicol}

\newtheorem{theorem}{Theorem}[section]
\newtheorem{prop}[theorem]{Proposition}
\newtheorem{lemma}[theorem]{Lemma}

\newtheorem{remark}{Remark}[section]

\newcommand\cA{{\cal A}}
\newcommand\cC{{\cal C}}

\newcommand\cI{{\cal I}}
\newcommand\cL{{\cal L}}
\newcommand\cB{{\cal B}}

\newcommand\cN{{\cal N}}

\newcommand\cP{{\cal P}}

\newcommand\cR{{\cal R}}

\newcommand\e{\epsilon}
\newcommand\ve{\varepsilon}

\newcommand\ov{\overline}

\def\bbr{{\mathbb R}}

\def\text#1{\hbox{#1}}

\def\endproof{\mbox{\ $\qed$}}

\def\E{{\bf E}}
\def\e{{\bf e}}
\def\A{{\bf A}}
\def\P{{\bf P}}
\def\C{{\bf C}}

\def\H{{\bf H}}
\def\L{{\bf L}}
\def\O{\hbox{\rm O}}

\def\Chi{{\bf 1}}
\def\d{\mathrm{d}}
\def\build #1_#2{\mathrel{\mathop{\kern 0pt #1}\limits_{#2}}} 
\newcommand{\zs}[1]{{\mathchoice{#1}{#1}{\lower.25ex\hbox{$\scriptstyle#1$}}
{\lower0.25ex\hbox{$\scriptscriptstyle#1$}}}}
\numberwithin{equation}{section}

\begin{document}

\title{Adaptive nonparametric estimation in heteroscedastic regression models.\\
{\Large Part 2: Asymptotic efficiency.}}
\author{{\Large By  Leonid Galtchouk and Sergey Pergamenshchikov}
\thanks{The second author is partially supported by the RFFI-Grant 04-01-00855.}\\
Louis Pasteur University  of Strasbourg and University of Rouen
}

\date{}
\maketitle

\begin{abstract}
In the paper we study asymptotic properties of the adaptive
procedure  proposed in the paper Galtchouk, Pergamenshchikov, 2007, 
for nonparametric estimation of unknown regression.
We prove that this
procedure is asymptotically efficient for some quadratic risk,
i.e. we show that the asymptotic quadratic risk for this procedure
coincides with the Pinsker constant which gives a sharp 
 lower bound for quadratic risk over all possible estimates.
\footnote{
{\sl AMS 2000 Subject Classification} : primary 62G08; secondary 62G05, 62G20}
\footnote{
{\sl Key words}: adaptive estimation, 
 asymptotic bounds, efficient estimation,
heteroscedastic regression, nonparametric estimation, 
non-asymptotic estimation, oracle inequality, Pinsker's constant.
}
\end{abstract}
\bibliographystyle{plain}
\renewcommand{\columnseprule}{.1pt}
\newpage

\section{Introduction}\label{I}

The paper deals with the  estimation problem in 
the heteroscedastic nonparametic regression model 
\begin{equation}\label{I.1}
y_j=S(x_j)+\sigma_j(S)\,\xi_j\,, 
\end{equation}
where the design points $x_j=j/n$, $S(\cdot)$
 is an unknown function to be estimated, 
$(\xi_j)_\zs{1\le j\le n}$ is a sequence of centered i.i.d. random variables with
unit variance and $\E\xi_1^4=\xi^*<\infty$,
$(\sigma_j(S))_\zs{1\le j\le n}$ are unknown scale functionals depending on
unknown regression function $S$ and the design points.


 Typically, the notion of asymptotic
optimality is associated with the optimal convergence rate of the minimax risk
(see for example, Ibragimov, Hasminskii,1981; Stone,1982).
  An important question in optimality
results is to study the exact asymptotic behaviour of the minimax risk. Such results
have been obtained only in a limited number of investigations. As to the
nonparametric estimation problem
for heteroscedastic regression models we should mention the papers Efromovich, 2007,
Efromovich, Pinsker, 1996, and Galtchouk, Pergamenshchikov, 2005,
 concerning the exact asymptotic behaviour of the $\cL_\zs{2}$-risk
and paper by Brua, 2007,
devoted to the efficient pointwise estimation for
heteroscedastic regressions. 

We remind that an example of 
heteroscedastic regression models is given by econometrics
(see, for example, Goldfeld, Quandt, 1972,
 p. 83), where 
for  consumer budget problems one  uses
 some parametric version of 
model \eqref{I.1}  with the scale coefficients  defined as
\begin{equation}\label{I.2}
\sigma^2_j(S)=c_\zs{0}+c_\zs{1}x_\zs{j}+c_\zs{2}S^2(x_\zs{j})\,,
\end{equation}
where $c_\zs{0}$, $c_\zs{1}$ and $c_\zs{2}$ are some positive unknown constants.

The purpose of the article is to study asymptotic properties of the adaptive 
estimation procedure proposed in Galtchouk, Pergamenshchikov, 2007,
for which a non-asymptotic oracle
inequality was proved for quadratic risks. We will prove that this oracle inequality
is  asymptotically sharp, i.e. the asymptotic quadratic risk is minimal. It means
the adaptive estimation procedure is efficient under some conditions on the scales
$(\sigma_j(S))_\zs{1\le j\le n}$ which are satisfied in the case \eqref{I.2}.
Note that in  Efromovich, 2007, Efromovich, Pinsker, 1996,
an efficient adaptive procedure is constructed
for heteroscedastic regression when the scale coefficient is independent of $S$, i.e.
$\sigma_j(S)=\sigma_j$. In Galtchouk, Pergamenshchikov, 2005,
for the model \eqref{I.1} the asymptotic efficiency
was proved under strong conditions on the scales which are not satisfied in the case \eqref{I.2}.
Moreover in the cited papers the efficiency was proved for the gaussian random variables
$(\xi_j)_\zs{1\le j\le n}$ that is very restrictive for applications of proposed methods
to practical problems.

In the paper we modify the risk by introducing into a additional supremum with respect to a classe
of unknown noise distributions like to Galtchouk, Pergamenshchikov, 2006.
 This modification allow us to eliminate from
the risk dependence on the noise distribution. Moreover for this risk a efficient procedure 
is robust with respect to changing of noise distributions. 

It is well known to prove the asymptotic efficiency one has to show that the asymptotic quadratic 
risk coincides with the lower bound which is equal to the Pinsker constant. In the paper two 
problems
are resolved: in the first one an upper bound for the risk is obtained by making use of the
non-asymptotic oracle inequality from Galtchouk, Pergamenshchikov, 2007,
in the second one we prove that this upper
bound coincides with the Pinsker constant. Let us remind that the adaptive procedure proposed in
Galtchouk, Pergamenshchikov, 2007, is based on weighted mean-squares estimates, where the 
weights are corresponding
modifications of the Pinsker weights for the homogene case (when
$\sigma_1(S)=\ldots=\sigma_n(S)=1$) relative to a certain smoothness of the function $S$ and
this procedure chooses an estimator best for the quadratic risk among these estimates. To obtain 
the Pinsker constant for the model \eqref{I.1} one has to prove a sharp asymptotic lower bound 
for the quadratic risk in the case when the noise variance depends on the unknown regression
function. This lower bound is obtained by making use of an inequality of kind of the van Trees
inequality (see, Gill, Levit, 1995).
 First we prove the inequality  for a parametric regression with 
the noise variance depending on the unknown regression (see Section 6) and further we apply the
inequality to the nonparametric regression by standard reducing to a parametric case.

    The paper is organized as follows.  
In Section~\ref{Ad} we construct a adaptive estimation procedure.
In Section~\ref{Co} we formulate principal conditions. The main result
is given in Section~\ref{M}. The upper bound for the quadratic risk is given
in Section~\ref{U}. In Section~\ref{Tr} we find the lower bound for a parametric model. 
In Section~\ref{Fa}
we study the parametric family. In Section~\ref{L} we obtain the lower bound  for
model \eqref{I.1}.
 An appendix contains some technical results.

\section{Adaptive procedure}\label{Ad}

In this section we describe the adaptive procedure proposed in 
\cite{GaPe1}. We make use of the 
 standard trigonometric basis $(\phi_j)_\zs{j\ge 1}$ in $\cL_2[0,1]$, i.e.
\begin{equation}\label{Ad.0}
\phi_1(x)=1\,,\quad
\phi_\zs{j}(x)=\sqrt{2}\,Tr_\zs{j}(2\pi [j/2]x)\,,\ j\ge 2\,,
\end{equation}
where the function $Tr_\zs{j}(x)=\cos(x)$ for even $j$ and
$Tr_\zs{j}(x)=\sin(x)$ for odd $j$; $[x]$ denotes the integer part of $x$.
We remind that if $n$ is odd then the functions $(\phi_j)_\zs{1\le j\le n}$ 
are orthonormal with respect to  the empirical inner product generated by the sieve
$(x_\zs{j})_\zs{1\le j\le n}$ in \eqref{I.1}, i.e. for any $1\le i,j\le n$,
$$
(\phi_i\,,\,\phi_j)_\zs{n}=
\frac{1}{n}\sum^n_\zs{l=1}\phi_i(x_l)\phi_j(x_l)={\bf Kr}_\zs{ij}\,,
$$
where ${\bf Kr}_\zs{ij}$ is Kronecker's symbol.
Thanks to this basis we pass to  the discrete Fourier transformation of model \eqref{I.1}, i.e.
\begin{equation}\label{Ad.1}
\hat{\vartheta}_\zs{j,n}=\vartheta_\zs{j,n}+(1/\sqrt{n})\xi_\zs{j,n}\,,
\end{equation}
where $\hat{\theta}_\zs{j,n}=(Y,\phi_j)_n$,
$\theta_\zs{j,n}=(S,\phi_j)_n$
and
$$
\xi_\zs{j,n}=\frac{1}{\sqrt{n}}\sum^n_\zs{l=1}\sigma_l(S)\xi_l\phi_j(x_l)\,.
$$
Here  $Y=(y_\zs{1},\ldots,y_\zs{n})'$ and
$S=(S(x_\zs{1}),\ldots,S(x_\zs{n}))'$. The prime denotes the transposition.

We estimate the function $S$ by the weighted least squares estimator
\begin{equation}\label{Ad.2}
\hat{S}_\zs{\lambda}=\sum^n_\zs{j=1}\lambda(j)\hat{\vartheta}_\zs{j,n}\phi_\zs{j}\,,
\end{equation}
where the weight vector $\lambda=(\lambda(1),\ldots,\lambda(n))'$
belongs to some finite set $\Lambda$ from $[0,1]^n$ with $n\ge 3$.
Here we make use of the weight family $\Lambda$ introduced in \cite{GaPe1}, i.e.

\begin{equation}\label{Ad.3}
\Lambda\,=\,\{\lambda_\zs{\alpha}\,,\,\alpha\in\cA_\zs{\ve}\}\,,
\quad
\cA_\zs{\ve}=\{1,\ldots,k_*\}\times\{t_1,\ldots,t_m\}\,,
\end{equation}
where $k_*=[1/\sqrt{\ve}]$, $t_i=i\ve$, $m=[1/\ve^2]$ and $\ve=1/\ln n$.

For any $\alpha=(\beta,t)\in\cA_\zs{\ve}$ we define the weight vector
$\lambda_\zs{\alpha}=(\lambda_\zs{\alpha}(1),\ldots,\lambda_\zs{\alpha}(n))'$ as
\begin{equation}\label{Ad.4}
\lambda_\zs{\alpha}(j)=\Chi_\zs{\{1\le j\le j_\zs{0}\}}+
\left(1-(j/\omega(\alpha))^\beta\right)\,
\Chi_\zs{\{ j_\zs{0}<j\le \omega(\alpha)\}}\,,
\end{equation}
where $j_0=j_\zs{0}(\alpha)=\left[\omega(\alpha)/\ln n\right]$,
$\omega(\alpha)=(A_\zs{\beta}\,t)^{1/(2\beta+1)}n^{1/(2\beta+1)}$ and
$$
A_\zs{\beta}=(\beta+1)(2\beta+1)/(\pi^{2\beta}\beta)\,.
$$
To find the optimal weights we choose the cost function equals to the penalized 
mean integrated squared error in which unknown parameters are replaced by some estimators.
The cost function is as follows
\begin{equation}\label{Ad.5}
J_\zs{n}(\lambda)\,=\,\sum^n_\zs{j=1}\,\lambda^2(j)\hat{\vartheta}^2_\zs{j,n}\,-
2\,\sum^n_\zs{j=1}\,\lambda(j)\,\tilde{\vartheta}_\zs{j,n}\,
+\,\rho \hat{P}_\zs{n}(\lambda)\,,
\end{equation}
where 
\begin{equation}\label{Ad.6}
\tilde{\vartheta}_\zs{j,n}=
\hat{\vartheta}^2_\zs{j,n}-\frac{1}{n}\hat{\varsigma}_\zs{n}
\quad\mbox{with}\quad
\hat{\varsigma}_\zs{n}=\sum^{n}_\zs{j=l_\zs{n}+1}
\hat{\vartheta}^2_\zs{j,n}
\end{equation}
and $l_\zs{n}=[n^{1/3}+1]$. The penalty term we define as
$$
\hat{P}_\zs{n}(\lambda)=\frac{|\lambda|^2 \hat{\varsigma}_\zs{n}}{n}\,,\quad
|\lambda|^2=\sum^n_\zs{j=1} \lambda^2(j)
\quad\mbox{and}\quad
\rho=\frac{1}{3+\ln^{\gamma} n}\,.
$$
for some $\gamma>0$.
Finally, we set
\begin{equation}\label{Ad.7}
\hat{\lambda}=\mbox{argmin}_\zs{\lambda\in\Lambda}\,J_n(\lambda)
\quad\mbox{and}\quad
\hat{S}_\zs{*}=\hat{S}_\zs{\hat{\lambda}}\,.
\end{equation}

The goal of this paper is to study asymptotic ($n\to\infty$) properties
of this estimation procedure.

\section{Conditions}\label{Co}

First we impose some conditions on unknown function $S$ in model \eqref{I.1}.

Let   $\cC^{k}_\zs{per,1}(\bbr)$ be  the set of $1$-periodic
$k$ times differentiable $\bbr\to\bbr $ functions. We assume that
 $S$ belongs to the following set
\begin{equation}\label{Co.1}
W^{k}_\zs{r}=\{f\in\cC^{k}_\zs{per,1}(\bbr)
\,:\,\sum_\zs{j=0}^k\,\|f^{(j)}\|^2\le r\}\,,
 \end{equation}
where $\|\cdot\|$ denotes the  norm in $\cL_\zs{2}[0,1]$, i.e.
\begin{equation}\label{Co.2}
\|f\|^2=\int^1_\zs{0}f^2(t)\d t\,.
\end{equation}
Moreover, we suppose that $r>0$ and $k\ge 1$ are unknown parameters.

Note that, we can represent the set $W^{k}_\zs{r}$ as
an ellipse in $\cL_\zs{2}[0,1]$, i.e.
\begin{equation}\label{Co.3}
W^{k}_\zs{r}=\{f\in\cL_\zs{2}[0,1]\,:\,
\sum_\zs{j=1}^\infty\,a_\zs{j}\vartheta^2_\zs{j}\le r\}\,,
 \end{equation}
where 
\begin{equation}\label{Co.3-1}
\vartheta_\zs{j}=(f,\phi_\zs{j})=\int^1_\zs{0}f(t)\phi_\zs{j}(t)\d t
\end{equation}
and 
\begin{equation}\label{Co.4}
a_\zs{j}=\sum^k_\zs{l=0}\|\phi^{(l)}_\zs{j}\|^2=
\sum^k_{i=0}(2\pi [j/2])^{2i}\,.
\end{equation}
Here $(\phi_\zs{j})_\zs{j\ge 1}$ is the trigonometric basis 
defined in \eqref{Ad.0}.

Now we decribe the conditions on the scale coefficients $(\sigma_j(S))_\zs{j\ge 1}$.

\begin{itemize}
\item[$\H_\zs{1})$] {\em 
 $\sigma_j(S)=g(x_j,S)$ for some unknown  function 
$g : [0,1]\times \L_1[0,1] \to \bbr_+$, which is
square integrable with respect to $x$ such that
\begin{equation}\label{Co.5}
\lim_\zs{n\to\infty}\,\sup_\zs{S\in W^k_r}\,\left|
n^{-1}\,\sum^n_\zs{j=1}\,g^2(x_j,S)\,-\,\varsigma(S)\,\right|\,=0\,,
\end{equation}
where $\varsigma(S):=\,\int_0^1\,g^2(x,S)\d x$. 
Moreover,
\begin{equation}\label{Co.6}
g_\zs{*}=\inf_\zs{0\le x\le 1}\,
\inf_\zs{S\in W^k_r} g^2(x,S)\,>0
\end{equation}
and
\begin{equation}\label{Co.6-1}
\sup_\zs{S\in W^k_r} \varsigma(S)<\infty\,.
\end{equation}
}
\item[$\H_\zs{2})$] {\em 
For any $x\in [0,1]$ the operator 
$g^2(x,\cdot)\,:\,\C[0,1]\to \bbr$
 is differentiable in  
the Fr\'echet sense for any fixed function $f_\zs{0}$ from $\C[0,1]$
, i.e.
 for any $f$ from some vicinity of $f_\zs{0}$ in $\C[0,1]$
$$
g^2(x,f)=g^2(x,f_\zs{0})+\L_\zs{x,f_\zs{0}}(f-f_\zs{0})+
\Upsilon(x,f_\zs{0},f)\,,
$$
where the Fr\'echet derivative
$\L_\zs{x,f_\zs{0}}\,:\,\C[0,1]\to \bbr$
is a bounded linear operator
and 
the residual term $\Upsilon(x,f_\zs{0},f)$ for each $x\in [0,1]$ satisfies the following
property
$$
\lim_\zs{|f-f_\zs{0}|_\zs{*}\to 0}
\frac{|\Upsilon(x,f_\zs{0},f)|}{|f-f_\zs{0}|_\zs{*}}=0\,,
$$
where $|f|_\zs{*}=\sup_\zs{0\le t\le 1} |f(t)|$.
}

\item[$\H_\zs{3})$] {\em
There exists some positive constant $C^*$ such that
for any function $S$ from $\C[0,1]$ the operator 
$\L_\zs{x,S}$ defined in condition $\H_\zs{2})$ 
satisfies the following inequality for any function $f$ from $\C[0,1]$
\begin{equation}\label{Co.8}
|\L_\zs{x,S}(f)|
\le C^*
\left(
|S(x)f(x)|+|f|_\zs{1}+\|S\|\,\|f\|
\right)\,,
\end{equation}
where $|f|_\zs{1}=\int^1_\zs{0}|f(t)|\d t$.
}

\item[$\H_\zs{4})$] {\em The function 
$g^2_\zs{0}(\cdot)=g^2(\cdot,S_\zs{0})$ corresponding to $S_\zs{0}\equiv 0$
is  continuous on the interval $[0,1]$.
Moreover,
$$
\lim_\zs{\delta\to 0}\,
\sup_\zs{0\le x\le 1}\,
\sup_\zs{|S|_\zs{*}\le \delta}\,
|g^2(x,S)-g^2(x,S_\zs{0})|\,=\,0\,.
$$
}
\end{itemize}
Now we give some examples of functions satisfying conditions $\H_\zs{1})$-$\H_\zs{4})$.

We fix some $c_\zs{0}>0$.  Let  $G\,:\,[0,1]\times\bbr\to [c_\zs{0}\,,\,+\infty)$
be a
 function such that
\begin{equation}\label{Co.9}
\lim_\zs{\delta\to 0}\max_\zs{|u-v|\le \delta}
\sup_\zs{y\in\bbr} 
| G(u,y)-G(v,y)|=0\,.
\end{equation}
and
\begin{equation}\label{Co.10}
G'_\zs{*}=
\sup_\zs{0\le x\le 1}
\sup_\zs{y\in\bbr}
|G_\zs{y}(x,y)|/|y|<\infty\,.
\end{equation}

Moreover,  let $V\,:\,\bbr\to\bbr_\zs{+}$ be a continuously differentiable  function such that
$$
v'_\zs{*}=
\,
\sup_\zs{y\in\bbr}\,
|\dot{V}(y)|/(1+|y|)<\infty\,.
$$
We set
\begin{equation}\label{Co.11}
g^2(x,S)=G(x,S(x))+\int^1_\zs{0}\,V(S(t))\d t\,.
\end{equation}
In this case 
$$
\varsigma(S)=\int^1_\zs{0}G(t,S(t))\d t+\int^1_\zs{0}\,V(S(t))\d t
$$
and for any $S\in W^k_\zs{r}$
\begin{align*}
\left| n^{-1}\sum^n_\zs{j=1}\,g^2(x_j,S)-\varsigma(S)\right|&
\le
\sum^n_\zs{j=1}\,
\int^{x_\zs{j}}_\zs{x_\zs{j-1}}
\left|G(x_j,S(x_\zs{j}))-
G(t,S(t))
\right|\,\d t\\
&\le
\Delta_\zs{n}+
\frac{G'_\zs{*}}{n}\,\int^1_\zs{0}\,|S(t)|\,|\dot{S}(t)|\d t\le \Delta_\zs{n}+\frac{G'_\zs{*}}{n}\, r\,,
\end{align*}
where
$$
\Delta_\zs{n}=\max_\zs{|u-v|\le 1/n}\sup_\zs{y\in\bbr}|G(u,y)-G(v,y)|\,.
$$
Therefore by condition \eqref{Co.9} we obtain $\H_\zs{1})$.

Moreover, the Fr\'echet derivative in this case is given  by
$$
\L_\zs{x,S}(f)=G_\zs{y}(x,S(x))f(x)+
\int^1_\zs{0}\dot{V}(S(t)) f(t)\d t\,.
$$
It is easy to see that this operator satisfies the inequality
  \eqref{Co.8} with 
$$
C^*=G'_\zs{*}+v'_\zs{*}\,.
$$
For example, we can take in \eqref{Co.11}
\begin{equation}\label{Co.12}
G(x,y)=c_\zs{0}+c_\zs{1}x+c_\zs{2}y^2
\quad\mbox{and}\quad
V(x)=c_\zs{3}x^2
\end{equation}
with some coefficients $c_\zs{0}>0$,  $c_\zs{i}\ge 0, i=1,2,3$.  
Therefore, we obtain the function \eqref{I.2} if we put in \eqref{Co.11}-\eqref{Co.12}
 $c_\zs{3}=0$, i.e. $V\equiv 0$.


\medskip
\vspace{5mm}
\section{Main results}\label{M}

Denote by $\cP_\zs{*}$  the family of unknown noise density. Remind that the noise
random variables $(\xi_\zs{j})_\zs{1\le j\le n}$ are centered with unit variance and
$\E\xi_\zs{1}^4\le \xi^*$, where $\xi^*\ge 3$.
For any estimate $\hat{S}$ we define the following quadratic risk
\begin{equation}\label{M.1}
\cR_n(\hat{S},S)=\sup_\zs{p\in\cP_\zs{*}}\E_\zs{S,p}\|\hat{S}-S\|^2_\zs{n}\,,
\end{equation}
where $\E_\zs{S,p}$ is the expectation with respect to the distribution $\P_\zs{S,p}$
of the observations $(y_\zs{1},\ldots,y_\zs{n})$ with the fixed function $S$ and
the fixed  density $p$ of random variables $(\xi_\zs{j})_\zs{1\le j\le n}$ in model \eqref{I.1},
$\|S\|^2_\zs{n}=(S,S)_\zs{n}$.

In Galtchouk, Pergamenshchikov, 2007,
we shown the following non-asymptotic Oracle inequality for procedure 
\eqref{Ad.7}.
\begin{theorem}\label{Th.M.1}
Assume that in  model \eqref{I.1} the function
$S$ belongs to  $W_\zs{r}^{1}$.
Then, for any odd $n\ge 3$, any $0<\rho<1/3$ and $r>0$,
 the estimate $\hat{S}_\zs{*}$ satisfies the following 
oracle inequality
\begin{equation}\label{M.2}
\cR_\zs{n}(\hat{S}_\zs{*},S)\,
\le(1+\kappa(\rho))
\min_\zs{\lambda\in\Lambda}\,
\cR_\zs{n}(\hat{S}_\zs{\lambda},S)\,
+n^{-1}\cB_\zs{n}(\rho)\,,
\end{equation}
where
$$
\kappa(\rho)=(6\rho-2\rho^2)/(1-3\rho)
$$
and the function $\cB_\zs{n}(\rho)$ is such that, for any $\delta>0$,
\begin{equation}\label{M.3}
\lim_\zs{n\to\infty}\cB_\zs{n}(\rho)/n^\delta=0\,.
\end{equation}
\end{theorem}

Now we formulate the main asymptotic results. To this end 
for any function  $S\in W^k_\zs{r}$ we set
\begin{equation}\label{M.4}
\gamma_k(S)\,=\,\Gamma^*_k\,r^{1/(2k+1)}\,(\varsigma(S))^{2k/(2k+1)}\,,
\end{equation}
where
$$
\Gamma^*_k=(2k+1)^{1/(2k+1)}\left(k/(\pi\,(k+1))\right)^{2k/(2k+1)}\,.
$$
It is well known (see, for example, Nussbaum, 1985) 
that for any function $S\in W^k_\zs{r}$
the optimal convergence rate is $n^{2k/(2k+1)}$.

\begin{theorem}\label{Th.M.2}
Assume that in model \eqref{I.1} the sequence
$(\sigma_j(S))$ fulfils the condition 
$\H_1)$. Then the estimator $\hat{S}_\zs{*}$ from \eqref{Ad.7}
 satisfies the inequality
\begin{equation}\label{M.5}
\limsup_{n\to\infty}\,
n^{2k/(2k+1)}\,\sup_\zs{S\in W^k_\zs{r}}\,
\cR_n(\hat{S}_\zs{*},S)/\gamma_k(S)\,\le 1\,.
\end{equation}
\end{theorem}
The following result  gives the sharp lower bound for risk \eqref{M.1} and show that
$\gamma_k(S)$ is \mbox{\it the Pinsker} \mbox{\rm constant}.
\begin{theorem}\label{Th.M.3}
Assume that in model \eqref{I.1}  the sequence
$(\sigma_j(S))$
satisfies the conditions $\H_2)$-- $\H_4)$. Then, for any estimate $\hat{S}_n$, the risk 
$\cR_n(\hat{S}_n,S)$
admits the following asymptotic lower bound
\begin{equation}\label{M.6}
\liminf_{n\to\infty}\,n^{2k/(2k+1)}\,\inf_{\hat{S}_n}\,\sup_\zs{S\in W^k_\zs{r}}\,
\cR_n(\hat{S}_n,S)/\gamma_k(S)\,
\ge 1\,.
\end{equation}
\end{theorem}

\begin{remark}\label{Re.M.1}
Note that in Galtchouk, Pergamenshchikov, 2005
an asymptotically efficient estimate was constructed and results
similar to Theorems \ref{Th.M.2} and \ref{Th.M.3} were claimed for the model  \eqref{I.1}.
In fact the upper bound is true there under some additional condition on the smoothness of
the function $S$, i.e. on the parameter $k$. In the cited paper this additional condition
is not formulated since erroneous inequality $(A.6)$. To avoid the use of this inequality
we modify the estimating procedure by introducing the penalty term 
$\rho\, \hat{P}_\zs{n}(\lambda)$ in the cost function \eqref{Ad.5}. By this way we remove
all additional conditions on the smoothness parameter $k$.
\end{remark}
\medskip
\section{Upper bound}\label{U}

In this section we prove Theorem~\ref{Th.M.2}.  To this end we will make use of oracle inequality
\eqref{M.2}. We have to find an estimator from the family
\eqref{Ad.2}-\eqref{Ad.3} for which we can show the upper bound \eqref{M.5}. We start with the construction
of  such an estimator. First we put
\begin{equation}\label{U.0}
\tilde{l}_\zs{n}=\inf\{i\ge 1\,:\,i\ve\ge \ov{r}(S)\}\wedge m
\quad\mbox{and}\quad
\ov{r}(S)=r/\varsigma(S)\,.
\end{equation}
Then we choose an index from the set $\cA_\ve$ as
$$
\tilde{\alpha}=(k,\tilde{t}_\zs{n})\,,
$$
where $k$ is the parameter of the set $W^k_\zs{r}$ and $\tilde{t}_\zs{n}=\tilde{l}_\zs{n}\ve$.
Finally, we set
\begin{equation}\label{U.0-1}
\tilde{S}=\hat{S}_\zs{\tilde{\lambda}}
\quad\mbox{and}\quad
\tilde{\lambda}=\lambda_\zs{\tilde{\alpha}}\,.
\end{equation}
Now we show the upper bound \eqref{M.5} for this estimator.
\begin{theorem}\label{Th.U.1} Assume that condition $\H_1)$ hold. Then 
\begin{equation}\label{U.1}
\limsup_{n\to\infty}\,n^{2k/(2k+1)}\,
\sup_{S\in W^k_r}\,
\cR_\zs{n}(\tilde{S},S)/\gamma_k(S)\,\le 1\,.
\end{equation}
\end{theorem}
\begin{remark}\label{Re.U.1}
Note that the estimator
$\tilde{S}$ belongs to estimate family \eqref{Ad.2}-\eqref{Ad.3}, but we can't use directly
 this estimator because the parameters $k$, $r$ and $\ov{r}(S)$ are unknown. We can use
this upper bound only through the oracle inequality \eqref{M.2} proved for 
procedure \eqref{Ad.7}.
\end{remark}


\noindent {\bf Proof.} 
To prove the theorem 
we will adapt to the  heteroscedastic case the corresponding proof
 from Nussbaum, 1985.

First, from \eqref{Ad.2}  we obtain that, for any $p\in\cP_\zs{*}$,
\begin{equation}\label{U.1-1}
\E_\zs{S,p}\,\|\tilde{S}-S\|_n^2=
\sum_{j=1}^{n}\,(1\,-\,\tilde{\lambda}_j)^2
\vartheta^2_\zs{j,n}
+
\frac{1}{n}
\sum_{j=1}^n\,\tilde{\lambda}_j^2\varsigma_\zs{j,n}\,,
\end{equation}
where
 $$
\varsigma_\zs{j,n}=\frac{1}{n}\,
\sum^n_\zs{l=1}\sigma^2_\zs{l}(S)\phi^2_\zs{j}(x_\zs{l})\,.
$$
Setting now
$\tilde{\omega}=\omega(\tilde{\alpha})$,
$\tilde{j}_0=[\tilde{\omega}/ \ln n]$,
$\tilde{j}_1=[\tilde{\omega} \ln n]$
and
$$
\varsigma_\zs{n}=\frac{1}{n}\,
\sum^n_\zs{l=1}\sigma^2_\zs{l}(S)\,,
$$
we rewrite \eqref{U.1-1} as follows
\begin{equation}\label{U.1-2}
\E_\zs{S,p}\,\|\tilde{S}-S\|_n^2=\sum_{j=\tilde{j}_0+1}^{\tilde{j}_1-1}
(1\,-\,\tilde{\lambda}_j)^2\vartheta^2_\zs{j,n}
+\varsigma_n\,n^{-1}
\sum_{j=1}^n\,\tilde{\lambda}_j^2+\Delta_1(n)+\Delta_2(n)
\end{equation}
with 
$$
\Delta_1(n)\,=\,\sum_{j=\tilde{j}_1}^{n}\,\vartheta^2_\zs{j,n}
\quad
\mbox{and}
\quad
\Delta_2(n)\,=\,n^{-1}\,\sum_{j=1}^n\,\tilde{\lambda}_j^2
\left(\varsigma_\zs{j,n}-\varsigma_\zs{n}\right)\,.
$$
Note that we have decomposed the first term in the right-hand of \eqref{U.1-1} into the sum
$$
\sum_{j=\tilde{j}_0+1}^{\tilde{j}_1-1}
(1\,-\,\tilde{\lambda}_j)^2\vartheta^2_\zs{j,n}\,+\,\Delta_1(n)\,.
$$
This decomposition allows us to show that $\Delta_1(n)$ is negligible and further to
approximate the first term by a similar term in which the coefficients $\vartheta_\zs{j,n}$
will be replaced by the Fourier coefficients $\vartheta_\zs{j}$ of the function $S$.

Taking into account the definition of $\omega(\alpha)$ in \eqref{Ad.4}
we can bound $\tilde{\omega}$ as
$$
\tilde{\omega} \ge (A_\zs{k})^{1/(2k+1)}\,
n^{1/(2k+1)}\,(\ln n)^{-1/(2k+1)}\,.
$$
Therefore, by
Lemma~\ref{Le.A.1} we obtain 
$$
\lim_{n\to\infty}\sup_{S\in W_r^k}\,n^{2k/(2k+1)}\,\Delta_1(n)=0\,.
$$
Let us consider now the next term $\Delta_2(n)$. We have
$$
|\Delta_2(n)|
=
\left|\frac{1}{n^2}\,\sum^n_{d=1}\,\sigma^2_d\,
\sum_{j=1}^n\,\tilde{\lambda}_j^2\,\ov{\phi}_\zs{j}(x_d)\right|
\le \frac{\sigma_\zs{*}}{n}
\sup_\zs{0\le x\le 1}
\left|\sum_{j=1}^n\,\tilde{\lambda}_j^2\,\ov{\phi}_\zs{j}(x)\right|\,,
$$
where $\ov{\phi}_\zs{j}(x)=\phi^2_\zs{j}(x)-1$. Now by Lemma~\ref{Le.A.2} and
definition \eqref{Ad.4}
we obtain directly the same property for $\Delta_2(n)$, i.e.
$$
\lim_{n\to\infty}\sup_{S\in W_r^k}\,n^{2k/(2k+1)}\,|\Delta_2(n)|=0\,.
$$
Setting 
$$
\hat{\gamma}_\zs{k,n}(S)=n^{2k/(2k+1)}
\sum_{j=\tilde{j}_0}^{\tilde{j}_1-1}(1-\tilde{\lambda}_j)^2\vartheta^2_\zs{j}
+
\varsigma_n n^{-1/(2k+1)}
\sum_{j=1}^n\,\tilde{\lambda}_j^2
$$
and applying the well-known inequality
$$
(a+b)^2\le (1+\delta)a^2+(1+1/\delta)b^2
$$
to the first term in the right-hand side of inequality \eqref{U.1-2} we obtain that, for any  $\delta>0$
and for any $p\in\cP_\zs{*}$,
\begin{align}\nonumber
\E_\zs{S,p}\,\|\tilde{S}-S\|_n^2 &\le(1+\delta)
\,\hat{\gamma}_\zs{k,n}(S)\,n^{-2k/(2k+1)}\\[2mm] \label{U.1-3}
&+\Delta_1(n)+
\Delta_2(n)+(1+1/\delta)\,
\Delta_3(n)\,,
\end{align}
where 
$$
\Delta_3(n)=\sum_{j=\tilde{j}_0+1}^{\tilde{j}_1-1}(\vartheta_\zs{j,n}\,-\,\vartheta_\zs{j})^2\,.
$$
Taking into account 
that $k\ge 1$ and that
$$
\tilde{j}_\zs{1}\le 
 (A_\zs{k})^{1/(2k+1)}\,
n^{1/(2k+1)} 
(\ln n)^{(2k+2)/(2k+1)}\,,
$$
we can show through Lemma~\ref{Le.A.3} 
that
$$
\lim_{n\to\infty}\sup_{S\in W_r^k}\,n^{2k/(2k+1)}\,\Delta_3(n)=0\,.
$$
Therefore inequality  \eqref{U.1-3}
 yields
$$
\limsup_{n\to\infty} n^{2k/(2k+1)}
\sup_{S\in W^k_r}
\cR_\zs{n}(\tilde{S},S)/\gamma_k(S)\le 
\limsup_{n\to\infty}
\sup_{S\in W^k_r}
\hat{\gamma}_\zs{k,n}(S)/\gamma_k(S)
$$
and to prove \eqref{U.1} it suffices to show that
\begin{equation}\label{U.3}
\limsup_{n\to\infty}
\sup_{S\in W^k_r}
\hat{\gamma}_\zs{k,n}(S)/\gamma_k(S)
\le 1\,.
\end{equation}
First it should be noted 
that definition \eqref{U.0}
and
inequalities \eqref{Co.6}-\eqref{Co.6-1} imply directly
$$
\lim_\zs{n\to\infty}\sup_{S\in W^k_r}
\left|
\tilde{t}_\zs{n}/\ov{r}(S)
-1
\right|=0\,.
$$
Moreover,
by the definition of $(\tilde{\lambda}_\zs{j})_\zs{1\le j\le n}$ 
for sufficiently large $n$ for which $\tilde{t}_\zs{n}\ge \ov{r}(S)$
we can calculate the following supremum
\begin{align*}
\sup_\zs{j\ge 1}
\,n^{2k/(2k+1)}
(1-\tilde{\lambda}_\zs{j})^2/(\pi j)^{2k}
&= \pi^{-2k}
(A_\zs{k}\tilde{t}_\zs{n})^{-2k/(2k+1)}\\
&\le \pi^{-2k}
(A_\zs{k}\ov{r}(S))^{-2k/(2k+1)}\,.
\end{align*}
Therefore, taking into account the definition of the coefficients $(a_\zs{j})_\zs{j \ge 1}$
in \eqref{Co.4} we obtain that
$$
\limsup_{n\to\infty}
n^{2k/(2k+1)}
\sup_{S\in W^k_r}
\sup_{j\ge \tilde{j}_0}\,\pi^{2k}
(A_\zs{k}\ov{r}(S))^{2k/(2k+1)}
(1-\tilde{\lambda}_j)^2/a_j
\le 1\,.
$$
Moreover, by definition \eqref{Ad.4} we get that
$$
\lim_{n\to\infty}
\sup_{S\in W^k_r}
\left|
\,
n^{-1/(2k+1)}\,\sum^n_\zs{j=1}\tilde{\lambda}^2_\zs{j}
-
(A_\zs{k}\ov{r}(S))^{1/(2k+1)}
\int^1_\zs{0}(1-z^k)^2\d z
\right|=0\,.
$$
Taking into account definition of $W^k_\zs{r}$ in \eqref{Co.3}
and condition \eqref{Co.5}
 we
obtain inequality \eqref{U.3}. Hence Theorem~\ref{Th.U.1}.
\endproof

Now Theorem~\ref{Th.M.1} and Theorem~\ref{Th.U.1} imply Theorem~\ref{Th.M.2}.
\medskip

\section{Lower bound for parametric heteroscedastic regression models} \label{Tr}

Let $(\bbr^n,\cB(\bbr^n),\P_\zs{\vartheta},\vartheta\in\Theta\subseteq\bbr^l)$ 
be a 
 statistical model relative to the observations $(y_j)_\zs{1\le j\le n}$ 
governed by the regression equation
\begin{equation}\label{Tr.1}
y_j\,=\,S_\zs{\vartheta}(x_j)\,+\,\sigma_j(\vartheta)\,\xi_j\,,
\end{equation}
where $\xi_1,\ldots,\xi_n$ are i.i.d. $\cN(0,1)$ random variables,
 $\vartheta=(\vartheta_1,\ldots,\vartheta_l)^\prime$ is a unknown parameter vector,
 $S_\zs{\vartheta}(x)$ is a unknown (or known) function 
and $\sigma_\zs{j}(\vartheta)=g(x_\zs{j},S_\vartheta)$, with the function
$g(x,S)$ defined in condition $\H_\zs{1})$. Assume that a prior distribution
$\mu_\zs{\vartheta}$ of the parameter $\vartheta$ in $\bbr^l$ is defined by the density
 $\Phi(\vartheta)$ of
the following form
$$
\Phi(\vartheta)\,=\,
\Phi(\vartheta_1,\ldots,\vartheta_l)=\prod_{i=1}^l\varphi_i(\vartheta_i)\,,
$$
where $\varphi_i$ is a continuously differentiable bounded density on $\bbr$ with
$$
\cI_i=\int_{\bbr}\frac{(\dot{\varphi}_i(z))^2}{\varphi_i(z)}\d z
<\infty\,.
$$
Let $\lambda(\cdot)$ be a continuously differentiable $\bbr^l\to \bbr$ function
such that, for any $1\le i\le l$,
\begin{equation}\label{Tr.1-1}
\lim_\zs{|\theta_\zs{i}|\to\infty}\,\lambda(\vartheta)\,\varphi_\zs{i}(\vartheta_\zs{i})=0
\quad\mbox{and}\quad
\int_{\bbr^l}\,\left|
\lambda^{\prime}_\zs{i}(\vartheta)
\right|\,
\Phi(\vartheta)\d \vartheta<\infty\,,
\end{equation}
where 
$$ 
\lambda^{\prime}_\zs{i}(\vartheta)=
(\partial/\partial\vartheta_i)\,\lambda(\vartheta)\,.
$$
Let $\hat{\lambda}_n$ be an estimator of $\lambda(\vartheta)$ based on
 observations $(y_j)_\zs{1\le j\le n}$.
 For any $\cB(\bbr^n\times\bbr^l)$ -
mesurable integrable function 
$G(x,\vartheta), x\in\bbr^n, \vartheta\in\bbr^l$, 
we set
$$
\tilde{\E}\,G(Y,\vartheta)\,=\,\int_{\bbr^l}\,\E_\zs{\vartheta}\,G(Y,\vartheta)\,
\Phi(\vartheta)\,\d \vartheta\,,
$$
where $\E_\zs{\vartheta}$ is the expectation with respect to the distribution 
$\P_\zs{\vartheta}$ of the vector $Y=(y_\zs{1},\ldots,y_\zs{n})$.
 Note that in this case
$$
\E_\zs{\vartheta}\,G(Y,\vartheta)\,=\,\int_\zs{\bbr^n}\,G(v,\vartheta)\,
f(v,\vartheta)\,\d v\,,
$$
where
\begin{equation}\label{Tr.2}
f(v,\vartheta)=\prod^n_\zs{j=1}\,\frac{1}{\sqrt{2\pi} \sigma_j(\vartheta)}
\exp\left\{-\,\frac{(v_j-S_\zs{\vartheta}(x_j))^2}
{2\sigma_j^2(\vartheta)}\right\}\,.
\end{equation}

We prove the following result.
\begin{theorem}\label{Th.Tr.1}
Assume that conditions $\H_1)-\H_2)$ hold.
Moreover, assume that 
 the function $S_\zs{\vartheta}(\cdot)$ is
uniformly over $0\le x\le 1$
 differentiable in $\cC[0,1]$ with respect to
$\vartheta_\zs{i},\ 1\le i\le l$, i.e. 
for any $1\le i\le l$
there exists
a function $S^{\prime}_\zs{\vartheta,i}\in \cC[0,1]$
such that
\begin{equation}\label{Tr.2-1}
\lim_\zs{h\to 0}\,
\max_\zs{0\le x\le 1}
\left|\left(
S_\zs{\vartheta+h\e_\zs{i}}(x)-S_\zs{\vartheta}(x)-S^{\prime}_\zs{\vartheta,i}(x)h\right)/h
\right|\,=0\,,
\end{equation}
where $\e_\zs{i}=(0,....,1,...,0)'$,  all coordinates are $0$, except the ith equals to $1$ .
Then for any square integrable estimator $\hat{\lambda}_n$ of 
$\lambda(\vartheta)$ and any $1\le i\le l$,
\begin{equation}\label{Tr.3}
\tilde{\E}(\hat{\lambda}_n-\lambda)^2\ge
\Lambda_i^2/(F_i\,+\,B_i\,+\,\cI_\zs{i})\,,
\end{equation}
where $\Lambda_i=
\int_{\bbr^l}\,
\lambda^{\prime}_\zs{i}(\vartheta)\,
\Phi(\vartheta)\d \vartheta$, $F_i=\sum_{j=1}^n
\int_\zs{\bbr^l}
\,(S^{\prime}_\zs{\vartheta,i}(x_j)/\sigma_j(\vartheta))^2\,
\Phi(\vartheta)\d \vartheta$
and
$$
 B_i=\,\frac{1}{2}\,\sum^n_\zs{j=1}\int_\zs{\bbr^l}
\frac{\tilde{\L}^2_i(x_j,S_\zs{\vartheta})}{\sigma^4_\zs{j}(S_\zs{\vartheta})}
\Phi(\vartheta)\d \vartheta\,, 
$$
$\tilde{\L}_i(x,\vartheta)=\L_\zs{x,S_\zs{\vartheta}}(S^{\prime}_\zs{\vartheta,i})$,
the operator $\L_\zs{x_,S}$ is defined in 
the condition $\H_2)$.
\end{theorem}
\noindent {\bf  Proof.} 
We put
$$
\varrho_i(v,\vartheta)=\frac{1}{f(v,\vartheta)\Phi(\vartheta)}
\,\frac{\partial}{\partial\vartheta_i}\left(f(v,\vartheta)\Phi(\vartheta)\right)\,.
$$
Note that due to condition \eqref{Co.6} the density \eqref{Tr.2} is bounded, i.e.
$$
f(v,\vartheta)\le (2\pi g_\zs{*})^{-n/2}\,.
$$
So through \eqref{Tr.1-1} we obtain that
$$
\lim_\zs{|\vartheta_\zs{i}|\to\infty}
\lambda(\vartheta)\,f(v,\vartheta)\varphi_\zs{i}(\vartheta_\zs{i})=0\,.
$$
Therefore, integrating by parts yields
\begin{align*}
\tilde{\E}(\hat{\lambda}_n-\lambda)\varrho_\zs{i}&=
\int_\zs{\bbr^{n+l}}
(\hat{\lambda}_n(v)-\lambda(\vartheta))
\frac{\partial}{\partial\vartheta_i}
\left(f(v,\vartheta)\Phi(\vartheta)\right)\d \vartheta\d v\\
&=
\int_\zs{\bbr^l}\left(\frac{\partial}{\partial\vartheta_i}\lambda(\vartheta)\right)
\Phi(\vartheta)\left(\int_\zs{\bbr^{n}}f(v,\vartheta)\d v \right)\d \vartheta=\Lambda_i\,.
\end{align*}
Now the Bouniakovskii-Cauchy-Schwarz inequality gives the following lower bound
$$
\tilde{\E}(\hat{\lambda}_n-\lambda)^2\ge
\Lambda_i^2/\tilde{\E}\varrho_\zs{i}^2\,.
$$
To estimate the denominator in the last ratio, note that
\begin{align*}
\varrho_\zs{i}(v,\vartheta)&=
\frac{1}{f(v,\vartheta)}
\frac{\partial}{\partial\vartheta_i} f(v,\vartheta)
+
\frac{\dot{\varphi}_\zs{i}(\vartheta_\zs{i})}{\varphi_\zs{i}(\vartheta_\zs{i})}\\
&=\tilde{f}_\zs{i}(v,\vartheta)
+
\frac{\dot{\varphi}_\zs{i}(\vartheta_\zs{i})}{\varphi_\zs{i}(\vartheta_\zs{i})}\,,
\end{align*}
where 
$$
\tilde{f}_\zs{i}(v,\vartheta)=
(\partial/\partial\vartheta_i)\ln f(v,\vartheta)\,.
$$
From \eqref{Tr.1} it follows that
$$
\tilde{f}_\zs{i}(v,\vartheta)=
\sum^{n}_\zs{j=1}(\xi^2_\zs{j}-1)\,
\frac{1}{2\sigma^2_j(\vartheta)}
\frac{\partial}{\partial\vartheta_i}\,\sigma^2_\zs{j}(\vartheta)
+
\sum^{n}_\zs{j=1}\,\xi_\zs{j}\,
\frac{S^{\prime}_\zs{i}(x_\zs{j})}{\sigma_j(\vartheta)}\,.
$$

Moreover, conditions $\H_\zs{2})$ and \eqref{Tr.2-1} imply
\begin{align*}
(\partial/\partial\vartheta_i)\,\sigma^2_\zs{j}(\vartheta)\,=\,
\partial/\partial\vartheta_i)\,g^2(x_\zs{j},S_\zs{\vartheta})\,=\,
\tilde{\L}_\zs{i}(x_j,
\vartheta)
\end{align*}
from which it follows
$$
\tilde{\E}\,\left(
\tilde{f}_\zs{i}(Y,\vartheta)
\right)^2\,
=F_\zs{i}+B_\zs{i}\,. 
$$
This implies inequality \eqref{Tr.3}. Hence Theorem~\ref{Th.Tr.1}.
\endproof

\medskip

\section{ Parametric kernel function family}\label{Fa}

In this section we define and study some special parametric kernel functions family 
which will be used to prove the sharp lower bound \eqref{M.6}.

Let us begin by kernel functions. We fix $\eta >0$ and we set 
\begin{equation}\label{Fa.1}
I_{\eta}(x)
=\eta^{-1}\int_{\bbr}\Chi_\zs{(|u|\le 1-\eta)}\,V\left(\frac{u-x}{\eta}\right)\,\d u\,,
\end{equation}
where $\Chi_A$ is the indicator of a set $A$, the kernel $V\in\C^{\infty}(\bbr)$  is 
such that
$$
V(u)=0
\quad\mbox{for}\quad |u|\ge 1
\quad\mbox{and}\quad 
\int^1_{-1}\,V(u)\,\d u=1\,.
$$
It is easy to see that the function $I_{\eta}(x)$ possesses the properties :
\begin{align*}
&0\le I_{\eta}\le 1\,,\quad I_{\eta}(x)=1
\quad\mbox{for}\quad |x|\le 1-2\eta \quad\mbox{and}\\ 
&I_{\eta}(x)=0
\quad\mbox{for}\quad |x|\ge 1\,.
\end{align*}
 Moreover,
for any $c>0$ and $m\ge 1$
\begin{equation}\label{Fa.2}
\lim_\zs{\eta\to 0}\,\sup_\zs{f\,:\,|f|_\zs{*}\le c}\,
\left|
\int_{\bbr}f(x)I^m_{\eta}(x)\d x-\int_{-1}^{1}f(x) \d x
\right|
=0\,,
\end{equation}
where $|f|_\zs{*}=\sup_\zs{-1\le x\le 1}|f(x)|$.

We divide the interval $[0,1]$ into $M$ equal parts of length $2h$ and on each of them
we construct a  kernel-type function that was used in Ibragimov, Hasminskii, 1981,
to obtain the lower bound
for estimation at a fixed point. A such constructed on each interval function equals to zero
at the extremities together with all derivatives. It means that Fourier partial sums with respect
to the trigonometric basis in $\cL_\zs{2}[-1,1]$ give a natural parametric approximation to the
function on each interval. 

Let  $(e_\zs{j})_\zs{j\ge 1}$ be the trigonometric basis in $L_2[-1,1]$, i.e.
\begin{equation}\label{Fa.3}
e_\zs{1}=1/\sqrt{2}\,,\quad
e_\zs{j}(x)=\,Tr_\zs{j}\left(\pi [j/2] x\right)\,,\ j\ge 2\,,
\end{equation}
where  $Tr_\zs{j}(x)=\cos(x)$ for even $j$ and
$Tr_\zs{j}(x)=\sin(x)$ for odd $j$. 

Now,
 for any array $z=\{(z_\zs{m,j})_\zs{1\le m\le M_\zs{n}\,,\,1\le j\le N_\zs{n}}\}$  we define
the following function
\begin{equation}\label{Fa.4}
S_\zs{z,n}(x)=\sum_{m=1}^{M_\zs{n}}\sum_{j=1}^{N_\zs{n}}\,z_\zs{m,j}\,D_\zs{m,j}(x)\,,
\end{equation}
where $D_\zs{m,j}(x)=e_j\left(v_m(x)\right)I_{\eta}\left(v_m(x)\right)$,
$$
v_m(x)=(x-\tilde{x}_m)/h_\zs{n}\,,
\quad\tilde{x}_m= 2mh_\zs{n}
\quad\mbox{and}\quad
M_\zs{n}=\left[1/(2h_\zs{n})\right]-1\,.
$$

We assume that the sequences 
$(N_\zs{n})_\zs{n\ge 1}$ and $(h_\zs{n})_\zs{n\ge 1}$,
satisfy the following conditions.

$\A_\zs{1})$
{\em The sequence $N_\zs{n}\to\infty$ as $n\to\infty$ and for any $p>0$
$$
\lim_\zs{n\to\infty}\,N^{p}_\zs{n}/n\,=\,0\,.
$$
Moreover, there exist $0<\delta_\zs{1}<1$ and $\delta_\zs{2}>0$
such that 
$$
h_\zs{n}=\,\O(n^{-\delta_\zs{1}})
\quad\mbox{and}\quad
h^{-1}_\zs{n}=\,\O(n^{\delta_\zs{2}})
\quad\mbox{as}\quad
n\to \infty\,.
$$
}
To define a prior distribution on the family of arrays,
we choose the following random array 
$\vartheta=\{(\vartheta_\zs{m,j})_\zs{ 1\le m\le M_\zs{n}\,,\, 1\le j\le N_\zs{n}}\}$ 
with
\begin{equation}\label{Fa.5}
\vartheta_\zs{m,j}\,=\,t_\zs{m,j}\,\zeta_\zs{m,j}\,,
\end{equation}
where  $(\zeta_\zs{m,j})$ are i.i.d. $\cN(0,1)$ random variables and 
$(t_\zs{m,j})_\zs{ 1\le m\le M_\zs{n}\,,\, 1\le j\le N_\zs{n}}$
are some nonrandom positive coefficients. We make use of gaussian variables since they
possess the minimal Fisher information and therefore maximize the lower bound
\eqref{Tr.3}.
 We set
\begin{equation}\label{Fa.6}
t^*_\zs{n}=\,\max_\zs{1\le m\le M_\zs{n}}
\sum^{N_\zs{n}}_\zs{j=1}\,t_\zs{m,j}\,.
\end{equation}
We assume that the coefficients $(t_\zs{m,j})_\zs{ 1\le m\le M_\zs{n}\,,\, 1\le j\le N_\zs{n}}$
satisfy the following conditions.

$\A_\zs{2})$
{\em There exists  a sequence of positive numbers $(d_\zs{n})_\zs{n\ge 1}$ such that
\begin{equation}\label{Fa.7}
\lim_\zs{n\to\infty}
\frac{d_\zs{n}}{h_\zs{n}^{2k-1}}\,
\sum^{M_\zs{n}}_\zs{m=1}\sum^{N_\zs{n}}_\zs{j=1}\,t^2_\zs{m,j}\,
j^{2(k-1)}=0\,,
\quad
\lim_\zs{n\to\infty}\,\sqrt{d_\zs{n}}\,t^*_\zs{n}=0\,,
\end{equation}
moreover, for any $p> 0$,
$$
\lim_\zs{n\to\infty}\,n^{p}\,\exp\{{-d_\zs{n}/2}\}=0\,.
$$
}

$\A_\zs{3})$
{\em For some $0<\varepsilon<1$
$$
\limsup_\zs{n\to\infty}
\frac{1}{h_\zs{n}^{2k-1}}\,\sum^{M_\zs{n}}_\zs{m=1}\sum^{N_\zs{n}}_\zs{j=1}\,t^2_\zs{m,j}\,
j^{2k}\,
\le (1-\varepsilon)r
\left(\frac{2}{\pi}\right)^{2k}
\,.
$$
}

$\A_\zs{4})$
{\em There exists $\epsilon_\zs{0}>0$ such that
$$
\lim_\zs{n\to\infty}
\frac{1}{h_\zs{n}^{4k-2+\epsilon_\zs{0}}}\,
\sum^{M_\zs{n}}_\zs{m=1}\sum^{N_\zs{n}}_\zs{j=1}\,t^4_\zs{m,j}\,j^{4k}\,=0\,.
$$
}

\begin{prop}\label{Pr.Fa.1}
Let conditions $\A_\zs{1})$--$\A_\zs{2})$. Then,
for any $p>0$ and for any $\delta>0$,
$$
\lim_\zs{n\to\infty}\,n^p\,\max_\zs{0\le l\le k-1}
\P\left(\|S^{(l)}_\zs{\vartheta,n}\|>\delta\right)=0\,.
$$
\end{prop}
\noindent{\bf Proof.} First note that for $0\le x\le 1$ we can represent the $l$th derivative as
\begin{equation}\label{Fa.8}
S_\zs{\vartheta,n}^{(l)}(x)=
\frac{1}{h^l}\sum_{m=1}^{M_\zs{n}}\,\sum^l_\zs{i=0}\,
\left(^{l}_{i}\right)\,I^{(l-i)}_\zs{\eta}(v_\zs{m}(x))
\,Q_\zs{i,m}(v_\zs{m}(x))\,,
\end{equation}
where 
$$
Q_\zs{i,m}(v)=
\sum_{j=1}^{N_\zs{n}}\vartheta_{m,j}e^{(i)}_j(v)\,.
$$
Therefore 
\begin{align*}
\|S^{(l)}_\zs{\vartheta,n}\|^2=
\frac{1}{h_\zs{n}^{2l-1}}\sum^{M_n}_\zs{m=1}
\int^1_\zs{-1}\,
\left(
\sum^l_\zs{i=0}\,
\left(^{l}_{i}\right)\,I^{(l-i)}_\zs{\eta}(v)
\,Q_\zs{i,m}(v)
\right)^2\,\d v
\end{align*}
and by the Bounyakovskii-Cauchy-Schwarz inequality we obtain that
\begin{equation}\label{Fa.9}
\|S^{(l)}_\zs{\vartheta,n}\|^2\le
\frac{C^*(l,\eta)}{h_\zs{n}^{2l-1}}
\sum^l_\zs{i=0}\,
\ov{Q}_\zs{i,m}
\end{equation}
with $C^*(l,\eta)=\max_\zs{-1\le v\le 1}\,\sum^l_\zs{i=0}\,
\left(\left(^{l}_{i}\right)\,I^{(l-i)}_\zs{\eta}(v)\right)^2$
and
$$
\ov{Q}_\zs{i,m}=\sum^{M_n}_\zs{m=1}
\int^1_\zs{-1}
\,Q^2_\zs{i,m}(v)
\,\d v
\,.
$$
Now we show that for any $0\le i \le k-1$ and $\delta>0$
\begin{equation}\label{Fa.10}
\lim_\zs{n\to\infty}
n^p\,
\P\left(
\ov{Q}_\zs{i,m}>\delta h_\zs{n}^{2k-1}
\right)=0\,.
\end{equation}
To that end  we introduce the following set
\begin{equation}\label{Fa.11}
\Xi_\zs{n}=\{\max_\zs{1\le m\le M_n}\,\max_\zs{1\le j\le N}\,\zeta^2_\zs{m,j}\le d_\zs{n}\}\,,
\end{equation}
where the sequence $(d_\zs{n})_\zs{n\ge 1}$ is given in condition $\A_\zs{2})$.
Therefore, taking into account that
\begin{align*}
\int^1_\zs{-1}
\,Q^2_\zs{i,m}(v)
\,\d v&=\sum^{N_\zs{n}}_\zs{j=1}\vartheta^2_\zs{m,j}
\int^1_\zs{-1}(e^{(i)}_\zs{j}(v))^2\,\d v\\
&\le \left(\frac{\pi}{2}\right)^{2i}
\sum^{N_\zs{n}}_\zs{j=1}t^2_\zs{m,j}\,j^{2i}\zeta^2_\zs{m,j}
\,,
\end{align*}
the function $\ov{Q}_\zs{i,m}$ can be estimated on the set 
$\Xi_\zs{n}$ as
$$
\ov{Q}_\zs{i,m}
\le\,
\left(\frac{\pi}{2}\right)^{2i}\,
d_\zs{n}
\sum^{M_n}_\zs{m=1}
\sum^{N_\zs{n}}_\zs{j=1}t^2_\zs{m,j}\,
j^{2i}\,
$$
and by  \eqref{Fa.7} we get, for any $\delta>0$ and
 sufficiently large $n$,
$$
\P\left(
\ov{Q}_\zs{i,m}
>\delta h_\zs{n}^{2k-1}\right)
\le
\,
\P\left(\Xi^c_\zs{n}\right)\,.
$$
Moreover, for sufficiently large $n$ 
$$
\P\left(\Xi^c_\zs{n}\right)\le \,M_\zs{n}\,N_\zs{n}\, e^{-d_\zs{n}/2}\,.
$$
Therefore,  conditions $\A_\zs{1})$ and  \eqref{Fa.7} imply 
\begin{equation}\label{Fa.12}
\limsup_\zs{n\to\infty}\,n^p\,
\P\left(\Xi^c_\zs{n}\right)=0\,,
\end{equation}
for any $p>0$. Hence Proposition~\ref{Pr.Fa.1}.

\endproof

\begin{prop}\label{Pr.Fa.2}
Let conditions
$\A_\zs{1})$--$\A_\zs{4})$. Then, for any $p>0$,
$$
\lim_\zs{n\to\infty}\,n^p\,
\P(S_\zs{\vartheta,n}\notin W^{k}_\zs{r})\,=0\,.
$$
\end{prop}
\noindent{\bf Proof.}
First of all we prove that for $\varepsilon$ from condition  $\A_\zs{3})$
\begin{equation}\label{Fa.13}
\lim_\zs{n\to\infty}\,n^p\,
\P\left(\|S^{(k)}_\zs{\vartheta,n}\|>\sqrt{(1-\varepsilon/4)r}\right)\,=0\,.
\end{equation}
Indeed, 
 putting in \eqref{Fa.8} $l=k$ we can represent the $k$th derivative of 
$S_\zs{\vartheta,n}$ as follows
\begin{equation}\label{Fa.14}
S_\zs{\vartheta,n}^{(k)}(x)=\hat{S}_\zs{k}(x)+
\ov{S}_\zs{k}(x)
\end{equation}
with
$$
\hat{S}_\zs{k}(x)=
\frac{1}{h^k}\sum_{m=1}^{M_\zs{n}}\,\sum^{k-1}_\zs{i=0}\,
\left(^{k}_{i}\right)\,I^{(k-i)}_\zs{\eta}(v_\zs{m}(x))
\,Q_\zs{i,m}(v_\zs{m}(x))
$$
and
$$
\ov{S}_\zs{k}(x)=
\frac{1}{h^k}\sum_{m=1}^{M_\zs{n}}\,I_\zs{\eta}(v_\zs{m}(x))
\,Q_\zs{k,m}(v_\zs{m}(x))\,.
$$
First, note that, we can estimate the norm of 
$\hat{S}_\zs{k}(x)$ by the same way as in inequality \eqref{Fa.9}, i.e.
$$
\|\hat{S}_\zs{k}\|^2
\le
\frac{C^*(k,\eta)}{h_\zs{n}^{2k-1}}
\sum^{k-1}_\zs{i=0}
\ov{Q}_\zs{i,m}\,.
$$
By making use of \eqref{Fa.10} we obtain 
that, for any $p>0$ and for any $\delta>0$,
\begin{equation}\label{Fa.15}
\lim_\zs{n\to\infty}\,n^p\,
\P\left(\|\hat{S}_\zs{k}\|>\delta\right)=0\,.
\end{equation}

Let us consider now the last term in \eqref{Fa.14}. Taking into account that
 $0\le I_\zs{\eta}(v)\le 1$ we get
\begin{align*}
\|\ov{S}_\zs{k}\|^2&=\frac{1}{h_\zs{n}^{2k-1}}
\sum^{M_\zs{n}}_\zs{m=1}
\int^1_\zs{-1}\,I^2_\zs{\eta}(v)Q^2_\zs{k,m}(v)\d v\\
&\le
\left(
\frac{\pi}{2}
\right)^{2k}
 \frac{1}{h_\zs{n}^{2k-1}}
\sum^{M_\zs{n}}_\zs{m=1}
\sum^{N_\zs{n}}_\zs{j=1}t^2_\zs{m,j}\,j^{2k}\,
\zeta^2_\zs{m,j}\,.
\end{align*}
Therefore from condition $\A_\zs{3})$ we get for sufficiently large $n$
$$
\|\ov{S}_\zs{k}\|^2\le (1-\varepsilon/2)r+\left(
\frac{\pi}{2}
\right)^{2k}
\sum^{M_\zs{n}}_\zs{m=1}\ov{\zeta}_\zs{m}:=
(1-\varepsilon/2)r+\left(\frac{\pi}{2}
\right)^{2k}Y_\zs{n}
$$
with
$$
\ov{\zeta}_\zs{m}=\frac{1}{h_\zs{n}^{2k-1}}
\sum^{N_\zs{n}}_\zs{j=1}t^2_\zs{m,j}\,
j^{2k}
\tilde{\zeta}_\zs{m,j}
\quad\mbox{and}\quad
\tilde{\zeta}_\zs{m,j}=\zeta^2_\zs{m,j}-1\,.
$$
We show that for any $p>0$ and for any $\delta>0$
\begin{equation}\label{Fa.16}
\lim_\zs{n\to\infty}\,n^p\,
\P\left(
\left|Y_\zs{n}\right|
>\delta\right)=0\,.
\end{equation}
Indeed, by the Chebyshev inequality for any $\iota>0$
\begin{equation}\label{Fa.17}
\P\left(
\left|Y_\zs{n}\right|
>\delta\right)\le\,
\E\,\left( Y_\zs{n}\right)^{2\iota}/\delta^{2\iota}\,.
\end{equation}
Note now that
 according to  the Burkholder-Davis-Gundy inequality for any $\iota>1$ there
exists a constant $B^*(\iota)>0$ such that
$$
\E\left( Y_\zs{n}\right)^{2\iota}\le B^*(\iota)\,
\E\,
\left(
\sum^{M_\zs{n}}_\zs{m=1}\ov{\zeta}^2_\zs{m}
\right)^{\iota}\,.
$$
Moreover, by putting
$$
\tilde{\zeta}_\zs{*}=\max_\zs{1\le m\le M_\zs{n}}\,
\max_\zs{1\le j\le N_\zs{n}}\,\tilde{\zeta}^2_\zs{m,j}
$$
we obtain that
$$
\ov{\zeta}^2_\zs{m}\,\le\,
\frac{N_\zs{n}}{h_\zs{n}^{4k-2}}
\sum^{N_\zs{n}}_\zs{j=1}t^4_\zs{m,j}\,j^{4k}\,
\tilde{\zeta}_\zs{*}\,.
$$
Therefore, by condition $\A_\zs{4})$ for sufficiently large $n$ 
\begin{align*}
\E\left( Y_\zs{n}\right)^{2\iota}&\le B^*(\iota)\,N^{\iota}_\zs{n}\,
h^{\iota\epsilon_\zs{0}}_\zs{n}\,\E\,\tilde{\zeta}^{\iota}_\zs{*}\\
&\le 
B^*(\iota)\,
\E\,(\zeta^2-1)^{2\iota}\,
M_\zs{n}N^{\iota+1}_\zs{n}
h^{\iota\epsilon_\zs{0}}_\zs{n}\,,
\end{align*}
where $\zeta\sim\cN(0,1)$.
Taking into account here condition $\A_\zs{1})$
 we obtain for sufficiently large $n$
$$
\E\left( Y_\zs{n}\right)^{2\iota}\le n^{-\delta_\zs{1}\,(\iota\epsilon_\zs{0}-2)}\,.
$$
Thus, choosing in  \eqref{Fa.17}
$$
\iota >p/(\epsilon_\zs{0}\delta_\zs{1})+2/\epsilon_\zs{0}
$$
 we obtain limiting equality 
\eqref{Fa.16} which together with \eqref{Fa.14}-\eqref{Fa.15}
implies \eqref{Fa.13}. Now it is easy to deduce that 
Proposition~\ref{Pr.Fa.1} yields Proposition~\ref{Pr.Fa.2}.
\endproof

\begin{prop}\label{Pr.Fa.3}
Let conditions
$\A_\zs{1})$--$\A_\zs{4})$. Then, for any $p>0$,
$$
\lim_\zs{n\to\infty}\,n^p\,
\E\,\|S_\zs{\vartheta,n}\|^2\,
\left(
\Chi_\zs{\{S_\zs{\vartheta,n}\notin W^{k}_\zs{r}\}}\,+\,
\Chi_\zs{\Xi^c_\zs{n}}
\right)
=0\,.
$$
\end{prop}
\noindent{\bf Proof.}
First of all, we remind that due to condition $\A_\zs{2})$
$$
\lim_\zs{n\to\infty}\,\sum^{M_\zs{n}}_\zs{m=1}\sum^{N_\zs{n}}_\zs{j=1}t^2_\zs{m,j}
\le\lim_\zs{n\to\infty}\,\frac{d_\zs{n}}{h^{2k-1}_\zs{n}}\,
\sum^{M_\zs{n}}_\zs{m=1}\sum^{N_\zs{n}}_\zs{j=1}t^2_\zs{m,j}\,j^{2(k-1)}=0\,.
$$
Therefore, taking into account that
\begin{equation}\label{Fa.17-1}
\|S_\zs{\vartheta,n}\|^2\le 
h_\zs{n}\sum^{M_\zs{n}}_\zs{m=1}\sum^{N_\zs{n}}_\zs{j=1}t^2_\zs{m,j}\zeta^2_\zs{m,j}
\end{equation}
we obtain, for sufficiently large $n$, 
$$
\E\,\|S_\zs{\vartheta,n}\|^2\,
\left(
\Chi_\zs{\{S_\zs{\vartheta,n}\notin W^{k}_\zs{r}\}}
+
\Chi_\zs{\Xi^c_\zs{n}}
\right)
\,
\le 
\max_\zs{m,j}\,
\E\,\zeta^2_\zs{m,j}
\left(
\Chi_\zs{\{S_\zs{\vartheta,n}\notin W^{k}_\zs{r}\}}
+
\Chi_\zs{\Xi^c_\zs{n}}
\right)
\,.
$$
Moreover, for any $1\le m\le M_\zs{n}$ and $1\le j\le N_\zs{n}$,
we estimate the last term as
\begin{align*}
\E\,\zeta^2_\zs{m,j}
\left(
\Chi_\zs{\{S_\zs{\vartheta,n}\notin W^{k}_\zs{r}\}}
+
\Chi_\zs{\Xi^c_\zs{n}}
\right)
&\le\,n\,\P(S_\zs{\vartheta,n}\notin W^{k}_\zs{r})\\[2mm]
&+
n\,\P(\Xi^c_\zs{n})
+
2\E\,\zeta^2\,\Chi_\zs{\{\zeta^2\ge n\}}\,,
\end{align*}
where $\zeta\sim \cN(0,1)$. By applying now Proposition~\ref{Pr.Fa.2}
and limit 
\eqref{Fa.12}
we obtain Proposition~\ref{Pr.Fa.3}.
\endproof

\begin{prop}\label{Pr.Fa.4}
Let conditions
$\A_\zs{1})$--$\A_\zs{4})$. Then for any function $g$
satisfying conditions \eqref{Co.6} and $\H_\zs{4})$ 
$$
\lim_\zs{n\to\infty}\,\sup_\zs{0\le x\le 1}\,
\E\,\left|\,g^{-2}(x,S_\zs{\vartheta,n})-g^{-2}_\zs{0}(x)\right|=0\,.
$$
\end{prop}
\noindent{\bf Proof.}
First, note that on the set $\Xi$ the random function $S_\zs{\vartheta,n}$
is uniformly bounded, i.e.
\begin{equation}\label{Fa.18}
|S_\zs{\vartheta,n}|_\zs{*}=
\sup_\zs{0\le x\le 1}\,|S_\zs{\vartheta,n}(x)|\,
\le \sqrt{d_\zs{n}}\,t^*_\zs{n}\,,
\end{equation}
where the coefficient $t^*_\zs{n}$ is defined in \eqref{Fa.6}. Therefore by condition $\H_\zs{1})$
we obtain
$$
\E\,\left|g^{-2}(x,S_\zs{\vartheta,n})- g^{-2}_\zs{0}(x)\right|\le 
\max_\zs{|S|_\zs{*}\le \sqrt{d_\zs{n}}\,t^*_\zs{n}}\,| g^{-2}(x,S)-g^{-2}_\zs{0}(x)|
+(2/g_\zs{*})\,
\P\left(\Xi^c_\zs{n}\right)\,.
$$
Conditions $\A_\zs{2})$ and $\H_\zs{4})$ together with the limit relation 
\eqref{Fa.12} imply Proposition~\ref{Pr.Fa.4}.
\endproof

\medskip
\section{Lower bound}\label{L}

In this section we prove Theorem~\ref{Th.M.3}.
To that end we establish the following auxiliary result.

\begin{lemma}\label{Le.L.1}
For any $0<\delta<1$ and  any estimate $\hat{S}_n$ of $S\in W^k_\zs{r}$,
$$
\|\hat{S}_n-S\|_n^2\ge(1-\delta)\|T_\zs{n}(\hat{S})-S\|^2\,
-\,(\delta^{-1}\,-\,1)\,r/n^2\,,
$$
where 
$T_\zs{n}(\hat{S})(x)=\sum_{k=1}^n\,\hat{S}_n(x_k)\Chi_\zs{(x_{k-1},x_k]}(x)$.
\end{lemma}
\noindent Proof of this Lemma is given in Appendix~\ref{Su.A.1}.

This Lemma implies that to prove  \eqref{M.6}, it suffices
to show the same asymptotic inequality for the integral risk, i.e. 
\begin{equation}\label{L.1}
\liminf_{n\to\infty}\,\inf_{\hat{S}_n}\,n^{2k/(2k+1)}\,
\cR_0(\hat{S}_n)\,\ge\,1\,,
\end{equation}
where
$$
\cR_0(\hat{S}_n)\,=\,
\sup_\zs{S\in W_r^k}\,\E_\zs{S,q}\,\|\hat{S}_n-S\|^2/\gamma_k(S)\,,
$$
$q$ is the gaussian $(0,1)$ density of the noise $(\xi_\zs{j})$ and 
$\|S\|^2=\int^1_\zs{0}S^2(x)\d x$.

To show \eqref{L.1} we will make use of the sequence of random functions
$(S_\zs{\vartheta,n})_\zs{n\ge 1}$ defined in \eqref{Fa.4}-\eqref{Fa.5}
with the coefficients $(t_\zs{m,j})$ satisfying conditions $\A_\zs{1})$--$\A_\zs{4})$
which will be chosen later.

 For any estimator $\hat{S}_n$, we denote by $\hat{S}^0_n$ its projection onto $W_r^k$, i.e.
 $\hat{S}^0_n=\hbox{\rm Pr}_\zs{W_r^k}(\hat{S}_n)$.
Since $W^k_\zs{r}$ is a convex set, we get that
$$
\|\hat{S}_n-S\|^2\ge\|\hat{S}^0_n-S\|^2\,.
$$
Therefore, we can write that
$$
\cR_0(\hat{S}_n)
\ge\int_{\{z:S_\zs{z,n}\in W^k_\zs{r}\}\cap\Xi_\zs{n}}
\,\frac{\E_\zs{S_\zs{z,n},q}\|\hat{S}^0_n-S_\zs{z,n}\|^2}{\gamma_k(S_\zs{z,n})}\,\mu_{\vartheta}(\d z)\,.
$$
Here $\mu_\zs{\vartheta}$ denotes 
 the distribution of $\vartheta$ in $\bbr^l$ with $l= M_\zs{n}N_\zs{n}$.
We recall also that the set $\Xi_\zs{n}$ is defined in \eqref{Fa.11}. 
Moreover, taking into account here inequality
\eqref{Fa.18} we estimate the risk $\cR_0(\hat{S}_n)$ from below as
$$
\cR_0(\hat{S}_n)
\ge\frac{1}{\gamma^*_\zs{n}}\,
\int_{\{z:S_\zs{z,n}\in W^k_\zs{r}\}\cap\Xi_\zs{n}}\,
\E_\zs{S_\zs{z,n},q}\|\hat{S}^0_n-S_\zs{z,n}\|^2\,\mu_{\vartheta}(\d z)
$$
with 
\begin{equation}\label{L.2}
\gamma^*_\zs{n}=\sup_\zs{|S|_\zs{*}\le \sqrt{d_\zs{n}}t^*_\zs{n}}\,\gamma_k(S)\,.
\end{equation}
Let us introduce now the corresponding Bayes risk
\begin{equation}\label{L.3}
\tilde{\cR}_0(\hat{S}^0_n)=
\int_\zs{\bbr^l}\,\E_\zs{S_\zs{z,n},q}\|\hat{S}^0_n-S_\zs{z,n}\|^2\,\mu_{\vartheta}(\d z)\,.
\end{equation}
Now through this risk we  rewrite  the lower bound for $\cR_0(\hat{S}_n)$ as 
\begin{equation}\label{L.4}
\cR_0(\hat{S}^0_n)
\ge
\tilde{\cR}_0(\hat{S}^0_n)/\gamma^*_\zs{n}-2\,
\varpi_\zs{n}/\gamma^*_\zs{n}
\end{equation}
with
$$
\varpi_\zs{n}=\E
(\Chi_\zs{\{S_\zs{\vartheta,n}\notin W^k_\zs{r}\}}\,+\,
\Chi_\zs{\Xi^c_\zs{n}})
(r+\|S_\zs{\vartheta,n}\|^2)\,.
$$
First of all, we reduce the nonparametric problem to parametric one. For this we replace the functions
$\hat{S}^0_n$ and $S$ by their Fourier series with respect to the basis
$$
\tilde{e}_\zs{m,i}(x)=(1/\sqrt{h})\,e_i\left(v_m(x)\right)\,
\Chi_\zs{\left(|v_m(x)|\le 1\right)}\,.
$$
By making use of this basis
we can estimate the norm $\|\hat{S}^0_n-S_\zs{z,n}\|^2$ from below as
$$
\|\hat{S}^0_n-S_\zs{z,n}\|^2 \ge
 \sum_{m=1}^{M_\zs{n}}\sum_{j=1}^{N_\zs{n}}\,(\hat{\lambda}_\zs{m,j}\,-\,
\lambda_\zs{m,j}(z))^2\,,
$$
where
$$
\hat{\lambda}_\zs{m,j}=\int_0^1\,\hat{S}^0_n(x)\tilde{e}_\zs{m,j}(x)\d x
\quad\mbox{and}\quad 
\lambda_\zs{m,j}(z)=\int_0^1\,S_\zs{z,n}(x)\tilde{e}_\zs{m,j}(x)\,\d x\,.
$$
Moreover, from definition \eqref{Fa.4}  one gets
$$
\lambda_\zs{m,j}(z)
=\sqrt{h}\sum_{i=1}^{N_\zs{n}}\,z_\zs{m,i}\int_{-1}^1\,e_i(u)e_j(u)I_{\eta}(u)\,\d u\,.
$$
It is easy to see that the functions $\lambda_\zs{m,j}(\cdot)$
satisfy condition \eqref{Tr.1-1} for gaussian prior densities. In this case
(see the definition in \eqref{Tr.3}) we have
$$
\Lambda_\zs{m,j}=
(\partial/\partial z_\zs{m,j}) \lambda_\zs{m,j}(z)
=\sqrt{h} \ov{e}_\zs{j}(I_\zs{\eta})\,,
$$
where 
\begin{equation}\label{L.6}
\ov{e}_\zs{j}(f)=\int_{-1}^1\,e^2_\zs{j}(v)\,f(v)\,\d v\,.
\end{equation}
Now to obtain a lower bound for
the Bayes risk 
$\tilde{\cR}_0(\hat{S}^0_n)$
 we make use of Theorem~\ref{Th.Tr.1} which implies that
\begin{equation}\label{L.7}
\tilde{\cR}_0(\hat{S}^0_n)\,\ge
\sum_{m=1}^{M_\zs{n}}\sum_{j=1}^{N_\zs{n}}\,
\frac{h \ov{e}^2_\zs{j}(I_\zs{\eta})}{F_\zs{m,j}+B_\zs{m,j}\,  
+t^{-2}_\zs{m,j}}\,,
\end{equation}
where $F_\zs{m,j}=\sum_{i=1}^n\,D^2_\zs{m,j}(x_i)\,
\E\,g^{-2}(x_\zs{i},S_\zs{\vartheta,n})$ and
$$
B_\zs{m,j}=
\frac{1}{2}\sum_{i=1}^n\,
\E\,
\left(\frac{\tilde{\L}_\zs{m,j}(x_\zs{i},S_\zs{\vartheta,n})}{g^2(x_\zs{i},S_\zs{\vartheta,n})}\right)^2\,
$$
with
$\tilde{\L}_\zs{m,j}(x,S)=\L_\zs{x,S} \left(D_\zs{m,j}\right)$.
In the appendix we show that
\begin{equation}\label{L.8}
\lim_\zs{n\to\infty}\,
\sup_\zs{1\le m\le M_\zs{n}}\,\sup_\zs{1\le j\le N_\zs{n}}
\left|
F_\zs{m,j}/(n h)\,-\,\ov{e}_j(I^2_{\eta})\,
g^{-2}_\zs{0}(\tilde{x}_m)
\right|=0
\end{equation}
and
\begin{equation}\label{L.9}
\lim_\zs{n\to\infty}
\sup_\zs{1\le m\le M_\zs{n}}\,\sup_\zs{1\le j\le N_\zs{n}}
\left|
B_\zs{m,j}/(n h)\right|\,=\,0\,.
\end{equation}
This means that, for any $\nu>0$
and for sufficiently large $n$,
$$
\sup_\zs{1\le m\le M_\zs{n}}\,\sup_\zs{1\le j\le N_\zs{n}}
\frac{F_\zs{m,j}+B_\zs{m,j}+t^{-2}_\zs{m,j}}
{nh \ov{e}_j(I^2_{\eta}) g^{-2}_\zs{0}(\tilde{x}_m)
+t^{-2}_\zs{m,j}} \le 1+\nu\,.
$$
Therefore, if we  denote in \eqref{L.7} 
$$
\kappa^2_\zs{m,j}=\,nh\,g^{-2}_\zs{0}(\tilde{x}_m)\,t^2_\zs{m,j}
\quad\mbox{and}\quad
\tau_\zs{j}(\eta,y)
=\ov{e}^2_\zs{j}(I_\zs{\eta})\,y/(\ov{e}_j^2(I^2_{\eta})y+1)
$$
we obtain that, for sufficiently large $n$,
$$
n^{2k/(2k+1)}\tilde{\cR}_0(\hat{S}^0_n)\,\ge
\frac{1}{1+\nu}\,n^{-1/(2k+1)}
\,
\sum_{m=1}^{M_\zs{n}}\,g^2_\zs{0}(\tilde{x}_m)
\,
\sum_{j=1}^{N_\zs{n}}\,\tau_\zs{j}(\eta,\kappa^2_\zs{m,j})\,.
$$
In the appendix we show that
\begin{equation}\label{L.10}
\lim_\zs{\eta\to 0}\,\sup_\zs{N\ge 1}\sup_\zs{(y_\zs{1},\ldots,y_\zs{N})\in \bbr^{N}_\zs{+}}
\left|
\sum_{j=1}^{N}\,\tau_\zs{j}(\eta,y_\zs{j})/
\sum_{j=1}^{N}\,\ov{\tau}(y_\zs{j})
\,-\,1
\right|\,=\,0\,,
\end{equation}
where
$$
\ov{\tau}(y)=y/(y+1)\,.
$$
Therefore we can write that, for sufficiently large $n$,
\begin{equation}\label{L.11}
n^{\frac{2k}{2k+1}}\tilde{\cR}_0(\hat{S}^0_n)\,\ge
\frac{1-\nu}{1+\nu}\,n^{-\frac{1}{2k+1}}
\,
\sum_{m=1}^{M_\zs{n}}\,g^2_\zs{0}(\tilde{x}_m)
\,J_\zs{N_\zs{n}}(\kappa^2_\zs{m,1},\ldots,\kappa^2_\zs{m,N_\zs{n}})\,,
\end{equation}
where
$$
J_\zs{N}(y_\zs{1},\ldots,y_\zs{N})=
\sum_{j=1}^{N}\,\ov{\tau}(y_\zs{j})\,.
$$
Obviously, to obtain a "good" lower bound for the risk 
$\tilde{\cR}_0(\hat{S}^0_n)$
 one needs to maximize the right-hand side of inequality \eqref{L.11}. Hence we
choose the coefficients $(\kappa^2_\zs{m,j})$ by maximization of the function
$J_\zs{N}$, i.e.
$$
\max_\zs{y_\zs{1},\ldots,y_\zs{N}}\,J_\zs{N}(y_\zs{1},\ldots,y_\zs{N})
\quad\mbox{subject to}\quad
\sum^N_\zs{j=1}y_\zs{j}j^{2k}\le R\,.
$$
The parameter $R>0$ will be chosen later to satisfy condition $\A_\zs{3})$. 
By the Lagrange multipliers method it is easy to find that
the solution of this problem is 
\begin{equation}\label{L.12}
y^*_\zs{j}(R)=(R+\sum^N_\zs{j=1}j^{2k})\,j^{-k}/\sum^N_\zs{j=1}j^{k}-1
\quad\mbox{for}\quad 1\le j\le N\,.
\end{equation}
To obtain a positive solution in \eqref{L.12} we need to impose  the following condition
\begin{equation}\label{L.13}
R\ge \,N^{k}\,\sum^N_\zs{j=1}j^{k}-\sum^N_\zs{j=1}j^{2k}\,.
\end{equation}
Moreover, from condition $\A_\zs{3})$ we  obtain that
\begin{equation}\label{L.14}
R\le 2^{2k+1}(1-\ve)r\,
n\,h^{2k+1}_\zs{n}/(\pi^{2k}\hat{g}_\zs{0}):=R^*_\zs{n}\,,
\end{equation}
where
$$
\hat{g}_\zs{0}=2h_\zs{n}\,\sum_{m=1}^{M_\zs{n}}\,g^2_\zs{0}(\tilde{x}_m)\,.
$$
Note that by condition $\H_\zs{4})$ the function 
$g_\zs{0}(\cdot)=g(\cdot,S_\zs{0})$ is continuous on 
the interval $[0,1]$, therefore 
\begin{equation}\label{L.15}
\lim_\zs{n\to\infty}\hat{g}_\zs{0}=\int^1_\zs{0}g^2(x,S_\zs{0})\d x=\varsigma(S_\zs{0})
\end{equation}
with $S_\zs{0}\equiv 0$.

Now we have to choose the sequence $(h_\zs{n})$. Note that if we put in \eqref{Fa.5}
\begin{equation}\label{L.15-1}
t_\zs{m,j}=(g_\zs{0}(\tilde{x}_m)/\sqrt{nh_\zs{n}})\,\sqrt{y^*_\zs{j}(R)}
\quad\mbox{i.e.}\quad
\kappa²_\zs{m,j}\,=\,y^*_\zs{j}(R)\,,
\end{equation}
we can rewrite inequality \eqref{L.11} as
\begin{equation}\label{L.16}
n^{\frac{2k}{2k+1}}\tilde{\cR}_0(\hat{S}^0_n)\,\ge
\frac{(1-\nu)}{(1+\nu)}\,\frac{\hat{g}_\zs{0}J^*_\zs{N_\zs{n}}(R)}{2h_\zs{n} n^{\frac{1}{2k+1}}}\,,
\end{equation}
where
$$
J^*_\zs{N}(R)=
N-\left(\sum^{N}_\zs{j=1}j^{k}\right)^2/(R+\sum^{N}_\zs{j=1}j^{2k})\,.
$$
It is clear that
$$
k^2/(k+1)^2\le \liminf_\zs{N\to\infty}\inf_\zs{R>0}J^*_\zs{N}(R)/N
\le \limsup_\zs{N\to\infty}\sup_\zs{R>0}J^*_\zs{N}(R)/N\le 1\,.
$$
Therefore to obtain a positive finite asymptotic lower bound in \eqref{L.16}
we have to take the parameter $h_\zs{n}$ as
\begin{equation}\label{L.17}
h_\zs{n}=h_\zs{*}n^{-1/(2k+1)}N_\zs{n}
\end{equation}
with some positive  coefficient $h_\zs{*}$. Moreover,  conditions 
\eqref{L.13}-\eqref{L.14} imply that
$$
(1-\ve)r\,\frac{2^{2k+1}}{\pi^{2k}}
\,\frac{1}{\hat{g}_\zs{0}}\,h^{2k+1}_\zs{*}
\ge
\frac{1}{N^{k+1}_\zs{n}}\sum^{N_\zs{n}}_\zs{j=1}j^k-\frac{1}{N^{2k+1}_\zs{n}}
\sum^{N_\zs{n}}_\zs{j=1}j^{2k}\,.
$$
Taking here limit as $n\to\infty$ thanks to asymptotic equality 
\eqref{L.15}, we obtain the following condition on $h_\zs{*}$
\begin{equation}\label{L.18}
h_\zs{*}\ge (\upsilon^*_\zs{\ve})^{1/(2k+1)}\,,
\end{equation}
where
$$
\upsilon^*_\zs{\ve}=\,\frac{k}{c^*_\zs{\ve}(k+1)(2k+1)}
\quad\mbox{and}\quad
c^*_\zs{\ve}=\frac{2^{2k+1}(1-\ve)r}{\pi^{2k}\varsigma(S_\zs{0})}\,.
$$
To maximize the function $J^*_\zs{N_\zs{n}}(R)$
at the right-hand side of inequality \eqref{L.16} we take  $R=R^*_\zs{n}$ defined in \eqref{L.14}.
Therefore we obtain that
\begin{equation}\label{L.19}
\liminf_\zs{n\to\infty}\,\inf_\zs{\hat{S}^0_n}\,
n^{2k/(2k+1)}\tilde{\cR}_0(\hat{S}^0_n)\,\ge
(\varsigma(S_\zs{0})/2)\,F(h_\zs{*})\,,
\end{equation}
where
$$
F(x)=\frac{1}{x}-\frac{2k+1}{(k+1)^2(c^*_\zs{\ve}(2k+1)x^{2k+2}+x)}\,.
$$
Taking into account that
$$
F^{\prime}(x)=-
\frac{(c^*_\zs{\ve}(2k+1)(k+1)x^{2k+1}-k)^2}{(k+1)^2(c^*_\zs{\ve}(2k+1)x^{2k+2}+x)^2}\le 0
$$
we find that
$$
\max_\zs{h_\zs{*}\ge (\upsilon^*_\zs{\ve})^{1/(2k+1)}}F(h_\zs{*})=
F((\upsilon^*_\zs{\ve})^{1/(2k+1)})=
(k/(k+1))(\upsilon^*_\zs{\ve})^{-1/(2k+1)}\,.
$$
This means that to obtain in \eqref{L.19} the maximal lower bound we have to take
in \eqref{L.17}
\begin{equation}\label{L.20}
h_\zs{*}=(\upsilon^*_\zs{\ve})^{1/(2k+1)}\,.
\end{equation}
Therefore, inequality \eqref{L.19} implies
\begin{equation}\label{L.21}
\liminf_\zs{n\to\infty}\,\inf_\zs{\hat{S}^0_n}\,
n^{2k/(2k+1)}\tilde{\cR}_0(\hat{S}^0_n)
\ge
(1-\ve)^{1/(2k+1)}\,\gamma_k(S_\zs{0})\,,
\end{equation}
where the function $\gamma_k(S_\zs{0})$ is defined in \eqref{M.4} for $S_\zs{0}\equiv 0$.

Now to end the definition of the sequence of the random functions
$(S_\zs{\vartheta,n})$ defined by \eqref{Fa.4} and 
\eqref{Fa.5} we have to define the sequence 
$(N_\zs{n})$.
 We remind that we make use of the  sequence 
$(S_\zs{\vartheta,n})$ 
with the coefficients
$(t_\zs{m,j})$
constructed in \eqref{L.15-1} for $R=R^*_\zs{n}$ given in 
\eqref{L.14} and for the sequence $h_\zs{n}$ given by 
\eqref{L.17} and \eqref{L.20} for some fixed arbitrary $0<\ve<1$.

  We will  choose the sequence $(N_\zs{n})$ 
to satisfy conditions $\A_\zs{1})$--$\A_\zs{4})$. We can take, for example 
$N_\zs{n}=[\ln^{4} n]+1$. Then condition $\A_\zs{1})$ is trivial. Moreover,
taking into account  that in
this case
$$
R^*_\zs{n}=\frac{2^{2k+1}(1-\ve)r}{\pi^{2k} \hat{g}_\zs{0}}\upsilon^*_\zs{\ve} N^{2k+1}_\zs{n}
=
\frac{\varsigma(S_\zs{0})}{\hat{g}_\zs{0}}
\frac{k}{(k+1)(2k+1)}\,N^{2k+1}_\zs{n}
$$
we find thanks to convergence
\eqref{L.15}
$$
\lim_\zs{n\to\infty}
(R^*_\zs{n}+\sum^{N_\zs{n}}_\zs{j=1}j^{2k})/(N^{k}_\zs{n} \sum^{N_\zs{n}}_\zs{j=1}j^{k})=1\,.
$$
Therefore,  solution \eqref{L.12} for sufficiently large $n$ 
satifies the following inequality
$$
\max_\zs{1\le j\le N_\zs{n}}
y^*_\zs{j}(R^*_\zs{n})\,j^k
\le 2 N^k_\zs{n}\,.
$$
Now it is easy to check  conditions 
 $\A_\zs{2})$ with $d_\zs{n}=\sqrt{N_\zs{n}}$ and $\A_\zs{4})$ for arbitrary 
$0<\epsilon_\zs{0}<1$.
As to condition $\A_\zs{3})$, note that by definition of $t_\zs{m,j}$ 
in \eqref{L.15-1}
we have 
\begin{align*}
\frac{1}{h_\zs{n}^{2k-1}}\,\sum^{M_\zs{n}}_\zs{m=1}\sum^{N_\zs{n}}_\zs{j=1}\,t^2_\zs{m,j}\,j^{2k}
&=\frac{1}{2 nh_\zs{n}^{2k+1}}\,\hat{g}_\zs{0}\,
\sum^{N_\zs{n}}_\zs{j=1}\,y^*_\zs{j}(R^*_\zs{n})\,j^{2k}\\
&
=\frac{R^*_\zs{n} \hat{g}_\zs{0}}{N^{2k+1}_\zs{n} 2\upsilon^*_\zs{\ve}}
=
(1-\varepsilon)r
\left(\frac{2}{\pi}\right)^{2k}\,.
\end{align*}
Hence condition $\A_\zs{3})$.

Therefore Propositions~\ref{Pr.Fa.2}-\ref{Pr.Fa.3} and limit \eqref{Fa.12}
imply that for any $p>0$
$$
\lim_\zs{n\to\infty}\,n^{p}\,
\varpi_\zs{n}=0\,.
$$
Moreover, by condition $\H_\zs{4})$ the sequence $\gamma^*_\zs{n}$ goes to $\gamma_k(S_\zs{0})$
as $n\to\infty$. Therefore, from this, \eqref{L.21} and \eqref{L.4} we get for any $0<\ve<1$
$$
\liminf_{n\to\infty}\,\inf_{\hat{S}_n}\,n^{2k/(2k+1)}\,
\cR_0(\hat{S}_n)\,\ge\,(1-\ve)^{1/(2k+1)}\,.
$$
Limiting here $\ve\to 0$ implies inequality \eqref{L.1}. 
Hence Theorem~\ref{Th.M.3}.
\endproof


\renewcommand{\theequation}{A.\arabic{equation}}
\renewcommand{\thetheorem}{A.\arabic{theorem}}
\renewcommand{\thesubsection}{A.\arabic{subsection}}
\section{Appendix}\label{Se.A}
\setcounter{equation}{0}
\setcounter{theorem}{0}

\subsection{Properties of trigonometric basis}\label{Su.A.0}

\begin{lemma}\label{Le.A.1} For any function $S\in  W_r^k$,
\begin{equation}\label{A.0-3}
\sup_{n\ge 1}\sup_{1\le m\le n-1}\,m^{2k}\,
\left(\sum_{j=m+1}^{n}\,\theta_\zs{j,n}^2
\right)
\,\le\,\frac{4r}{\pi^{2(k-1)}}\,.
\end{equation}
\end{lemma}

\begin{lemma}\label{Le.A.2}
For any $m\ge 0$,
\begin{equation}\label{A.0-2}
\sup_{N\ge 2}\quad\sup_{x\in [0,1]}N^{-m}\left|\sum_{l=2}^{N}\,l^m
\ov{\phi}_\zs{l}(x)\right|\le 2^m\,,
\end{equation}
where $\ov{\phi}_\zs{l}(x)=\phi^2_\zs{l}(x)-1$.
\end{lemma}
Proofs of Lemma~\ref{Le.A.1} and Lemma~\ref{Le.A.2} are given in \cite{GaPe1}.

\begin{lemma}\label{Le.A.3}
Let $\theta_\zs{j,n}$ and $\theta_j$ be the Fourier coefficients defined in \eqref{Ad.1}
and \eqref{Co.3-1} respectively.
 Then, for $1\le j\le n$
and $n\ge 2$,
\begin{equation}\label{A.0-1}
\sup_{S\in W^1_r}\,|\theta_\zs{j,n}-\theta_\zs{j}|\,\le\, 2\pi\,\sqrt{r}\,j/n\,.
\end{equation}
\end{lemma}
\noindent{\bf Proof.}
Indeed, we have
\begin{align*}
|\theta_\zs{j,n}-\theta_\zs{j}|&\,=\,
\left|\sum_{l=1}^n\int_{x_\zs{l-1}}^{x_\zs{l}}\,\left(S(x_l)\phi_j(x_l)-S(x)\phi_j(x)\right)\d x\right|\\
&\le\,n^{-1}\,\sum_{l=1}^n\int_{x_\zs{l-1}}^{x_\zs{l}}\,
\left(|\dot{S}(z)\phi_j(z)|\,+ \,|S(z)\dot{\phi}_j(z)|\right)\d z
\\
&=\,n^{-1}\,\int_{0}^{1}\,
\left(|\dot{S}(z)|\,|\phi_j(z)|\,+ \,|S(z)|\,|\dot{\phi}_j(z)|\right)\d z\,.
\end{align*}
By making use of the Bounyakovskii-Cauchy-Schwarz inequality we get
\begin{align*}
|\theta_\zs{j,n}-\theta_\zs{j}|&\le\,
n^{-1}\,\left(\|\dot{S}\|\,\|\phi\|\,+\,\|\dot{\phi}\|\,\|S\|\right)\\
&\le\,
n^{-1}\,\left(\|\dot{S}\|\,+\,\pi\,j\,\,\|S\|\right)\,.
\end{align*}
The definition of class $W^1_r$  implies (\ref{A.0-1}).
Hence Lemma~\ref{Le.A.1}.
\endproof

\subsection{Proof of Lemma~\ref{Le.L.1}}\label{Su.A.1}

First notice that, for any $S\in W_r^k$, one has
\begin{align*}
\|\hat{S}_n-S\|_n^2
=\|T_\zs{n}(\hat{S})-S\|^2+\Psi_n+\Delta_n\,,
\end{align*}
where
$$
\Psi_n=2\sum_{j=1}^n\int_\zs{x_{j-1}}^{x_j}(\hat{S}_n(x_j)-S(x))(S(x)-S(x_j))\d x
$$ 
and
$$
\Delta_n=\sum_{j=1}^n\int_\zs{x_{j-1}}^{x_j}\,(S(x)-S(x_j))^2\d x\,.
$$

For any $0<\delta<1$, by  making use of the elementary inequality
 $$
2ab\le \delta a^2+\delta^{-1}b^2\,,
$$
 one gets
$$
\Psi_n\,\le\, \delta\|T_\zs{n}(\hat{S})-S\|^2+\delta^{-1}\Delta_n\,.
$$
Moreover,  for any $S\in W_r^k$ with $k\ge 1$,
by the Bounyakovskii-Cauchy-Schwarz inequality we obtain that
\begin{align*}
\Delta_n\le\frac{1}{n}\sum_{j=1}^n\int_\zs{x_{j-1}}^{x_j}\,\dot{S}^2(t)\,\d t
=\frac{1}{n^2}\|\dot{S}\|^2\le\frac{r}{n^2}\,.
\end{align*}
Hence Lemma \ref{Le.L.1}. 
\endproof

\subsection{Proof of \eqref{L.8}}\label{Su.A.3}

First of all, note that Proposition~\ref{Pr.Fa.4},
condition \eqref{Co.6} and condition $\H_\zs{4})$ imply that
\begin{equation}\label{A.3-1}
\lim_\zs{n\to\infty}
\max_\zs{1\le m\le M_\zs{n}}\sup_\zs{0\le x\le 1}
\Chi_\zs{\{|v_m(x)|\le 1\}}
\,\E\left|g^{-2}(x,S_\zs{\vartheta,n})\,-\,g^{-2}_\zs{0}(\tilde{x}_m)
\right|\,=\,0\,.
\end{equation}
Let us show now that for any  continuously differentiable function
$f$ on $[-1,1]$
\begin{equation}\label{A.3-2}
\lim_\zs{n\to\infty}\,\sup_\zs{1\le m\le M_\zs{n}}\,
\left|
\frac{1}{nh}\,\sum^{n}_\zs{i=1}\,f(v_m(x_\zs{i}))\Chi_\zs{\{|v_m(x_\zs{i})|\le 1\}}
-
\int^{1}_\zs{-1}f(v)\d v
\right|=\,0\,.
\end{equation}
Indeed, setting 
$$
\Delta_\zs{n,m}=\frac{1}{nh}\,\sum^{n}_\zs{i=1}\,f(v_m(x_\zs{i}))\Chi_\zs{\{|v_m(x_\zs{i})|\le 1\}}
-
\int^{1}_\zs{-1}f(v)\d v
$$
 we obtain that
\begin{align*}
\left|
\Delta_\zs{n,m}
\right|&=
\left|
\frac{1}{nh}\,\sum^{i^*}_\zs{i=i_\zs{*}}\,f(v_m(x_\zs{i}))
-
\int^{1}_\zs{-1}f(v)\d v
\right|\\
&\le
\sum^{i^*}_\zs{i=i_\zs{*}}\,\int^{v_m(x_\zs{i})}_\zs{v_m(x_\zs{i-1})}|f(v_m(x_\zs{i}))-f(z)|\d z
+\max_\zs{|z|\le 1}|f(z)|(2-v^*+v_\zs{*})\,.
\end{align*}
where $i_\zs{*}=[n\tilde{x}_\zs{m}-nh]+1$, $i_\zs{*}=[n\tilde{x}_\zs{m}+nh]$, 
$$
v_\zs{*}=( [n\tilde{x}_\zs{m}-nh]+1-n\tilde{x}_\zs{m})/(nh)
\quad
\mbox{and}
\quad
v^*=([n\tilde{x}_\zs{m}+nh]-n\tilde{x}_\zs{m})/(nh)\,.
$$ 
Therefore, taking into accout that the derivative of the function $f$ is bounded
on the interval $[-1,1]$ we obtain that
$$
\left|
\Delta_\zs{n,m}
\right|\le
3\max_\zs{|z|\le 1}|\dot{f}(z)|/(nh_\zs{n})+
2\max_\zs{|z|\le 1}|f(z)|/(nh_\zs{n})\,.
$$
Taking into account the conditions on the sequence $(h_\zs{n})_\zs{n\ge 1}$
given in $\A_\zs{1})$ we obtain limiting equality \eqref{A.3-2} which together
with \eqref{A.3-1} implies \eqref{L.8}.
\endproof

\subsection{Proof of \eqref{L.9}}\label{Su.A.4}

Now we study the behaviour of $B_\zs{m,j}$. Due to inequality 
\eqref{Co.8}
we obtain that 
\begin{align*}
|\tilde{\L}_\zs{m,j}(x,S_\zs{\vartheta,n})|\le\,C^*\,
\left(|S_\zs{\vartheta,n}(x)D_\zs{m,j}(x)|+|D_\zs{m,j}|_\zs{1}
+
\|S_\zs{\vartheta,n}\|\,\|D_\zs{m,j}\|
\right)\,.
\end{align*}
Note that
\begin{align*}
\E(S_\zs{\vartheta,n}(x)D_\zs{m,j}(x))^2&=
\E\left(\sum^{N_\zs{n}}_\zs{l=1}\vartheta_\zs{m,l}e_\zs{l}(v_\zs{m}(x))\right)^2
e^2_\zs{j}(v_\zs{m}(x))I^4_\zs{\eta}(v_\zs{m}(x))\\
&\le 
\sum^{N_\zs{n}}_\zs{l=1}\,t^2_\zs{m,l}\,\Chi_\zs{\{|v_\zs{m}(x)|\le 1\}}
\le (t^*_\zs{n})^2\Chi_\zs{\{|v_\zs{m}(x)|\le 1\}}\,.
\end{align*}
We remind that the sequence $t^*_\zs{n}$ is defined in \eqref{Fa.6}.
Therefore, property \eqref{A.3-2} implies 
$$
\max_\zs{1\le m\le M_\zs{n}}\,\max_\zs{1\le j\le N_\zs{n}}\,
\frac{1}{nh}\,\sum^{n}_\zs{i=1}\,
\E(S_\zs{\vartheta,n}(x_\zs{i})D_\zs{m,j}(x_\zs{i}))^2
=\O((t^*_\zs{n})^2)\,.
$$
Moreover, as to the function $D_\zs{m,j}(\cdot)$ we find that
$$
|D_\zs{m,j}|_\zs{1}=\int^1_\zs{0}|e_\zs{j}(v_\zs{m}(x))\,
I_\zs{\eta}(v_\zs{m}(x))|\d x
=h\int^1_\zs{-1}| e_\zs{j}(v)\,
I_\zs{\eta}(v)|\d v\le 2h\,.
$$
Similarly we obtain $\|D_\zs{m,j}\|^2\le h$. 

Finally, by \eqref{Fa.17-1}we obtain that 
$$
\E\|S_\zs{\vartheta,n})\|^2\le h\sum^{M_\zs{n}}_\zs{m=1}
\sum^{N_\zs{n}}_\zs{j=1}\,t^2_\zs{m,j}\le (t^*_\zs{n})^2\,.
$$
Therefore,
$$
\max_\zs{1\le m\le M_\zs{n}}\,\max_\zs{1\le j\le N_\zs{n}}
B_\zs{m,j}/(nh)\,=\,\O((t^*_\zs{n})^2+h_\zs{n})
$$
and condition $\A_\zs{1})$ implies \eqref{L.9}.

\endproof
\subsection{Proof of \eqref{L.10}}\label{Su.A.5}

Indeed, by the direct calculation it easy to see that for any $N\ge 1$ and for any vector
$(y_\zs{1},\ldots,y_\zs{N})'\in\bbr^N_\zs{+}$
\begin{align*}
\left|
\frac{\sum_{j=1}^{N}\,\tau_\zs{j}(\eta,y_\zs{j})}
{\sum_{j=1}^{N}\,\ov{\tau}(y_\zs{j})}
\,-\,1
\right|\,\le\,
\frac{
\max_\zs{j\ge 1}\left(
|\ov{e}^2_\zs{j}(I_\zs{\eta})-\ov{e}_\zs{j}(I^2_\zs{\eta})|
+|\ov{e}^2_\zs{j}(I_\zs{\eta})-1|
\right)
}{\min_\zs{j\ge 1}\,\ov{e}_\zs{j}(I^2_\zs{\eta})}\,,
\end{align*}
where the operator $\ov{e}_\zs{j}(f)$ is defined in
in \eqref{L.6}.
Moreover, we remind that $\int^1_\zs{-1}e^2_\zs{j}(v)\d v=1$.
Therefore, taking into account  property
\eqref{Fa.2} we obtain \eqref{L.10}.
\endproof


\begin{flushright}
\begin{tabular}{lcl}
   L.Galtchouk                       &$\quad$& S. Pergamenshchikov              \\
 Department of Mathematics           &$\quad$& Laboratoire de Math\'ematiques Raphael Salem,\\      
 Strasbourg University               &$\quad$& Avenue de l'Universit\'e, BP. 12,            \\
 7, rue Rene Descartes               &$\quad$&  Universit\'e de Rouen,                  \\
 67084, Strasbourg, France           &$\quad$&  F76801, Saint Etienne du Rouvray, Cedex France.\\
 e-mail: galtchou@math.u-strasbg.fr  &$\quad$& Serge.Pergamenchtchikov@univ-rouen.fr         \\
\end{tabular}
\end{flushright} 

\end{document}